\documentclass[10pt,reqno,final]{amsart}
\usepackage{bm}
\usepackage{amsmath,amssymb,amsthm}

\newcommand{\D}{\mathbb{D}}

\newcommand{\x}{\mathbf{x}}

\newcommand{\ip}[2]{\left(#1,#2\right)}

\newcommand{\nff}{\mathbf{n}_f}
\newcommand{\npp}{\mathbf{n}_p}
\newcommand{\Ptau}{P_\tau}
\DeclareMathOperator{\tr}{tr}

\newcommand{\Omegaf}{\Omega_f}

\newcommand{\norm}[1]{\lVert #1 \rVert}
\newcommand{\xivec}{\boldsymbol{\xi}}
\newcommand{\etavec}{\boldsymbol{\eta}}

\usepackage{bm}
\newcommand{\Omf}{\Omega_f}
\newcommand{\Omp}{\Omega_p}

\newcommand{\ddt}{\frac{d}{dt}}

\usepackage{hyperref}
\usepackage{cleveref}
\usepackage{mathrsfs}
\usepackage{pifont}
\usepackage[english]{babel}
\usepackage{graphicx}
\usepackage{subfig}
\usepackage{hyperref}
\usepackage{enumitem}
\usepackage[notcite,notref]{showkeys}
\usepackage{color}
\usepackage{xcolor}
\usepackage{tabularx}
\usepackage{booktabs}
\usepackage{array}
\usepackage{makecell}
\usepackage{multirow,multicol}
\usepackage{caption}
\usepackage{algorithm,algorithmicx}
\usepackage{algpseudocode}
\usepackage{microtype}
\usepackage{fancybox}
\usepackage[foot]{amsaddr}
\usepackage{cite}
\allowdisplaybreaks

\usepackage{geometry}

\theoremstyle{plain}
\newtheorem{theorem}{Theorem}[section]
\newtheorem{lemma}{Lemma}[section]

\theoremstyle{definition}

\theoremstyle{remark}
\newtheorem{remark}{Remark}[section]

\numberwithin{equation}{section}
\graphicspath{{./figures/}}

\title[The time-dependent thermal-fluid-solid multiphysical coupled model]{Well-posedness of weak solutions for the time-dependent thermal-fluid-solid multiphysical coupled model}

\author[Yi Qin]{Yi Qin$^{1}$}
\author[Siyu Yang]{Siyu Yang$^{1}$}
\author[Luling Cao]{Luling Cao$^{1}$}
\author[Jian Li]{Jian Li$^{1,*}$}
\address{$^{1}$School of Mathematics and Date Science, Shaanxi University of Science and Technology, Xian,Shaanxi
710021, China}
\address{$^{1,*}$School of Mathematics and Data Science, Shaanxi University of science and technology, Xian, Shaanxi 710021,
China}
\email{yiqin@sust.edu.cn; 15594140979@163.com; lulingcao@sust.edu.cn; jianli@sust.edu.cn}

\thanks{This work is supported by the Innovative Team Project of Shaanxi Provincial Department of Education (25JP022); the Shaanxi Provincial Natural Science Basic Research Program Project (2026JC-YBMS-0006) and the Key Scientific Research Program of Shaanxi Provincial Department of Education (25JS019);  NSFC (12601779) and the Key Laboratory Project of the Shaanxi Provincial Department of Education (25JS020); and NSFC (12431016, 12471408) and the Shaanxi Provincial Science and Technology Innovation Team Project (2026RS-CXTD-71).}

\date{\today}

\begin{document}

\begin{abstract}
This paper constructs a time-dependent thermal-fluid-solid multiphysical coupled model describing the interaction between a free fluid and a thermo-poroelastic medium and provides a rigorous mathematical analysis. The model consists of the incompressible Navier-Stokes equations and the heat equation in the free-fluid region, together with the thermo-poroelasticity equations in the porous medium. The two regions are coupled on a fixed interface through conditions for normal velocity, tangential velocity, stress balance, temperature continuity, and heat flux continuity. Under a small-data assumption, we prove the existence and
uniqueness of a weak solution for this fully dynamic nonlinear coupled system. The proof employs the
Galerkin method to construct semi-discrete approximations, establishes uniform a priori estimates that
enable the extension of local solutions to global ones, and invokes compactness arguments to pass to the
limit, yielding a solution to the original coupled system.

\end{abstract}
\maketitle
\section{Introduction}

Fluid-poroelastic interaction is encountered in a wide range of applications, including geothermal systems \cite{shanker2022convective}, nuclear waste disposal \cite{grambow2014nuclear,liu2014non}, permafrost engineering \cite{deng2019application}, and biomechanics \cite{guo2016biomechanical}. Although temperature effects within the porous medium have been extensively incorporated into thermo-poroelasticity models, the thermal coupling between the free-fluid region and the porous matrix remains largely unaddressed in theoretical analysis. In geothermal extraction, for example, the heat exchange between the injected cold fluid and the high-temperature reservoir directly influences production efficiency, yet existing theoretical models cannot fully describe this mechanism. Motivated by this consideration, the present paper constructs a time-dependent thermal-fluid-solid multiphysical coupled model that fully characterizes the interaction among free-fluid flow, heat transfer, and thermo-poroelastic deformation across a fixed interface, and provides a rigorous well-posedness analysis of its weak solutions.

Several mathematical and numerical studies have addressed fluid-poroelastic interaction problems \cite{kunwar2020second,wang2022locking,guo2022decoupled}. For the Navier-Stokes-Biot coupling, Li, Caucao, and Yotov \cite{li2024augmented} introduced a partially augmented fully-mixed formulation and proved existence, uniqueness, stability, and error estimates; Wang and Yotov \cite{wang2025lagrange} developed a model with a Lagrange multiplier and proved existence, uniqueness, and stability under a small-data assumption; Cesmelioglu \cite{cesmelioglu2017analysis} considered a weak formulation of the time-dependent incompressible Navier-Stokes and Biot equations and proved existence, uniqueness, and a priori estimates under a similar small-data assumption; Wen et al. \cite{wen2021discontinuous} developed semi-discrete and fully discrete schemes with corresponding error estimates. However, thermal effects are not considered in these studies. In contrast, thermo-poroelasticity models have witnessed significant advances on multiple fronts, ranging from theoretical analysis to numerical methods. Santos, Carcione, and Ba \cite{santos2021existence} introduced coupling temperature terms and relaxation terms to describe wave propagation, proving the existence, uniqueness, and regularity of solutions in bounded domains; Jin \cite{jin2018thermo} revealed the significant effects of temperature gradients on infiltration kinetics and liquid content variation for the processing of interpenetrating phase composites. Nordbotten and Keilegavlen \cite{nordbotten2021introduction} developed a unified MPxA framework applicable to general polyhedral grids; Yi and Lee \cite{yi2024physics} developed a fully coupled enriched Galerkin scheme that achieves local conservation for mechanics, flow, and heat sub-problems at lower cost; Chen, Cui, and Zhou \cite{chen2025symmetric} proposed a symmetric interior penalty discontinuous Galerkin method for nonlinear quasi-static problems, including proofs of existence, uniqueness, and a priori error estimates; Cai et al. \cite{cai2026parameter} developed four-field formulations and robust preconditioners to ensure efficient solution of the discrete system under mesh refinement and parameter changes. Nevertheless, existing thermo-poroelasticity models focus solely on the porous medium and do not account for coupling with an adjacent free fluid.

In this work, we establish a time-dependent thermal-fluid-solid multiphysical coupled model that describes the interaction between an incompressible free fluid and a thermoporoelastic medium: the free-fluid region is governed by the Navier-Stokes equations and the heat equation, while the porous region is described by a thermo-poroelasticity model accounting for solid skeleton deformation, pore pressure evolution, and heat transfer; the two regions are coupled across a fixed interface through continuity conditions for velocity, stress, temperature, and heat flux, including the Beavers-Joseph-Saffman condition for tangential slip. The main contribution of this paper is to provide a rigorous mathematical analysis of this fully coupled nonlinear system: we establish the existence and uniqueness of a weak solution under a suitable smallness assumption, employing the Galerkin method combined with compactness arguments and a priori estimates, and recover the fluid pressure via an inf-sup condition.

The remainder of this paper is structured as follows. \Cref{sec:2} presents the full governing equations.
\Cref{sec:3} establishes the existence of a weak solution using a Galerkin method, while \cref{sec:4} proves
its uniqueness. \Cref{sec:5} reconstructs the fluid pressure via the inf-sup condition. The main theorem
is summarized in \cref{sec:6}, and conclusions are drawn in \cref{sec:7}.

\section{The time-dependent thermal-fluid-solid multiphysical coupled model}\label{sec:2}
\begin{figure}[htbp]
    \centering
    \includegraphics[width=0.5\linewidth, trim=3cm 11cm 5cm 10cm, clip]{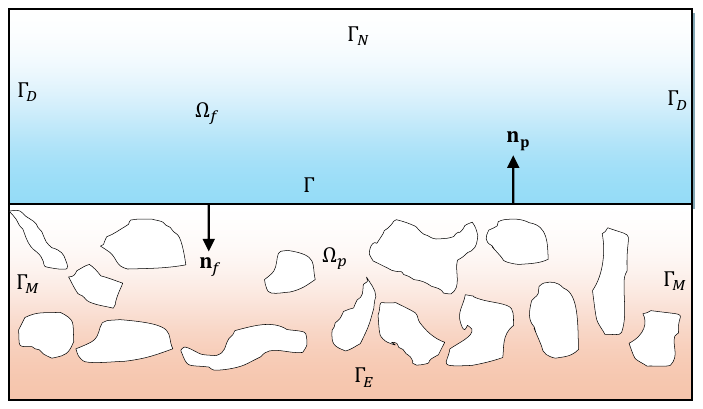}
        \caption{A global domain $\Omega$ consisting of a fluid flow region $\Omega_f$ and a porous media flow region $\Omega_p$ separated by an interface $\Gamma$.}
        \label{fig:1}
    \end{figure}
Consider a bounded domain $\Omega \subset \mathbb{R}^d$ $(d=2,3)$ consisting of a fluid region $\Omega_f$ and a porous media region $\Omega_p$, separated by an interface $\Gamma = \partial \Omega_f \cap \partial \Omega_p$ (see Figure~\ref{fig:1}). For a finite time interval $[0,T]$, the governing equations are described as follows.
The fluid motion with heat transfer is governed by the time-dependent incompressible Navier-Stokes equations and the heat equation\cite{abdullah2020coupled,boland1990analysis,burns2016feedback,shang2013error}. Find the fluid velocity $\mathbf{u}_f(\mathbf{x},t)$, pressure $p_f(\mathbf{x},t)$ and temperature $\theta_f(\mathbf{x},t)$ such that
\begin{align}
\partial_t \mathbf{u}_f + (\mathbf{u}_f\cdot\nabla)\mathbf{u}_f - \nabla\cdot \mathbb{T}_f(\mathbf{u}_f,p_f) &= \frac{Gr}{Re^2} \boldsymbol{\xi} \,\theta_f && \text{in }\Omega_f\times(0,T], \label{eq:fluid_momentum}\\
\nabla\cdot \mathbf{u}_f &= 0 && \text{in }\Omega_f\times(0,T], \label{eq:fluid_continuity}\\
\partial_t \theta_f - a_f \Delta \theta_f + (\mathbf{u}_f\cdot\nabla)\theta_f &= \mathbf{u}_f\cdot\boldsymbol{\xi} && \text{in }\Omega_f\times(0,T], \label{eq:fluid_heat}
\end{align}
with initial conditions
\begin{align}
\mathbf{u}_f(\mathbf{\x},0)=\mathbf{u}_f^0(\x),\qquad \theta_f(\x,0)=\theta_f^0(\x) \quad \text{in }\Omega_f,
\end{align}
and boundary conditions on $\partial\Omega_f\setminus\Gamma = \Gamma_{D}\cup\Gamma_{N}$:
\begin{equation}
\begin{aligned}
\mathbf{u}_f &= 0 && \text{on } \partial\Omega_f\setminus\Gamma\times[0,T],\\
\theta_f &= 0 && \text{on } \Gamma_{N}\times[0,T],\qquad
\frac{\partial\theta_f}{\partial\mathbf{n}_f} = 0 && \text{on } \Gamma_{D}\times[0,T].
\end{aligned}
\end{equation}
Here, the Cauchy stress tensor is given by $\mathbb{T}_f(\mathbf{u}_f,p_f)=2\mu_f\mathbb{D}(\mathbf{u}_f)-p_f\mathbb{I}$, with $\mathbb{D}(\mathbf{u}_f)=\frac12(\nabla\mathbf{u}_f+(\nabla\mathbf{u}_f)^\mathrm{T})$ being the strain rate tensor. The dynamic viscosity is $\mu_f=\frac{1}{Re}$, while $a_f$ denotes the thermal conductivity, $Gr$ the Grashof number, and $\boldsymbol{\xi} = (0,\dots,0,-1) \in \mathbb{R}^d$, where the first $d-1$ components are zero, which is the unit vector aligned with gravity.

In the porous domain $\Omega_p$, the energy conservation equation, momentum balance equation, and mass conservation equation for the displacement $\boldsymbol{\eta}_p(\mathbf{x},t)$, velocity $\boldsymbol{\xi}_p(\mathbf{x},t)$, pressure $p_p(\mathbf{x},t)$, and temperature $\theta_p(\mathbf{x},t)$ are given by \cite{chen2022multiphysics,jin2025thermo}
\begin{align}
\partial_t\bigl(a_0\theta_p - b_0 p_p + \beta\nabla\cdot\etavec_p\bigr) - \nabla\cdot(\bm{a}_p\nabla\theta_p) &= \phi && \text{in }\Omega_p\times(0,T], \label{eq:porous_heat}\\
\partial_t\xivec_p - \nabla\cdot\mathbb{T}_p(\etavec_p,p_p) + \beta\nabla\theta_p &= \mathbf{f} && \text{in }\Omega_p\times(0,T], \label{eq:porous_momentum}\\
\partial_t\bigl(S_0 p_p - b_0\theta_p + \alpha\nabla\cdot\etavec_p\bigr) - \nabla\cdot(\mathbf{K}\nabla p_p) &= g && \text{in }\Omega_p\times(0,T], \label{eq:porous_pressure}
\end{align}
with initial conditions
\begin{align}
\etavec_p(\mathbf{x},0)=\etavec_p^0(\mathbf{x}),\quad \theta_p(\mathbf{x},0)=\theta_p^0(\x),\quad p_p(\mathbf{x},0)=p_p^0(\mathbf{x}) \quad \text{in }\Omega_p,
\end{align}
and boundary conditions on $\partial\Omega_p\setminus\Gamma = \Gamma_{M}\cup\Gamma_{E}$
\begin{equation}
\begin{aligned}
\etavec_p &= 0 && \text{on } \Gamma_{E}\times[0,T], \quad&
\bigl(\mathbb{T}_p(\etavec_p,p_p)+\beta\nabla\theta_p\bigr)\cdot\mathbf{n}_p = 0 && \text{on } \Gamma_{M}\times[0,T],\\
\theta_p &= 0 && \text{on } \Gamma_{E}\times[0,T], &
\frac{\partial\theta_p}{\partial\mathbf{n}_p} = 0 && \text{on } \Gamma_{M}\times[0,T],\\
p_p &= 0 && \text{on } \Gamma_{E}\times[0,T], &
\mathbf{K}\nabla p_p\cdot\mathbf{n}_p = 0 && \text{on } \Gamma_{M}\times[0,T].
\end{aligned}
\end{equation}
Here $\xivec_p = \partial_t \etavec_p$. The total stress tensor $\mathbb{T}_p$ is defined as
$\mathbb{T}_p(\etavec_p,p_p)=2\mu_p\mathbb{D}(\etavec_p)+\lambda_p(\nabla\cdot\etavec_p)\mathbb{I}-\alpha p_p\mathbb{I}-\beta \theta_p\mathbb{I}$,
and the strain tensor is $\mathbb{D}(\etavec_p)=\frac12(\nabla\etavec_p+(\nabla\etavec_p)^\mathrm{T})$.
Moreover, $a_0$ denotes the effective thermal capacity, $b_0$ the thermal dilation coefficient, $\beta$ the thermal stress coefficient, $\mathbf{K}$ the permeability tensor, $\bm{a}_p$ the effective thermal conductivity, $\mu_p$ and $\lambda_p$ the Lamé parameters, $\alpha$ the Biot-Willis constant, and $S_0$ the specific storage coefficient.

Finally, on $\Gamma \times (0,T]$ we impose the following interface conditions: continuity of normal flux which satisfies mass conservation, the Beavers–Joseph–Saffman condition, balance of normal fluid stress components, balance of stresses, and continuity of temperature and heat flux \cite{hansbo2005nitsche,burman2009stabilization,hirata2007stability}.
\begin{equation}
\begin{aligned}
&\mathbf{u}_f\cdot\mathbf{n}_f = (\partial_t\etavec_p - \mathbf{K}\nabla p_p)\cdot\mathbf{n}_f, \\
&P_{\boldsymbol{\tau}}\bigl(\mathbb{T}_f(\mathbf{u}_f,p_f)\cdot\mathbf{n}_f\bigr) = -\frac{\alpha\sqrt{d}}{\sqrt{\operatorname{tr}(\mathbf{K})}} P_{\boldsymbol{\tau}}(\mathbf{u}_f - \partial_t\etavec_p), \\
&\mathbf{n}_f\cdot\bigl(\mathbb{T}_f(\mathbf{u}_f,p_f)\mathbf{n}_f\bigr) = -p_p, \\
&\mathbb{T}_f(\mathbf{u}_f,p_f)\mathbf{n}_f = \mathbb{T}_p(\etavec_p,p_p)\mathbf{n}_f, \\
&\theta_f = \theta_p, \\
&\mathbf{n}_f\cdot a_f\nabla\theta_f = -\mathbf{n}_p\cdot \bm{a}_p\nabla\theta_p,\label{eq:interface_heatflux}
\end{aligned}
\end{equation}
where $\{\boldsymbol{\tau}_f^i\}_{i=1}^{d-1}$ is a set of linearly independent tangent vectors on $\Gamma$ \cite{guo2022decoupled}, and $P_{\boldsymbol{\tau}}(\mathbf{v})=\sum_{i=1}^{d-1}(\mathbf{v}\cdot\boldsymbol{\tau}_f^i)\boldsymbol{\tau}_f^i$ denotes the tangential projection.

Now, we introduce the standard Lebesgue and Sobolev spaces \cite{qin2022variable,cesmelioglu2017analysis}:
\begin{align*}
\mathbf{H}_f &= \{\mathbf{v}_f\in H^1(\Omega_f)^2 : \mathbf{v}_f=0 \text{ on } \partial\Omega_f\},\\
X_f &= L^2_0(\Omega_f)= \{q \in L^2(\Omega_f) : \int_{\Omega_f} q \,\mathrm{d}x = 0\},\\
T_f &= \{\varphi\in H^1(\Omega_f) : \varphi=0 \text{ on } \Gamma_{D}\},\\
\mathbf{H}_p &= \{\mathbf{v}_p\in H^1(\Omega_p)^2 : \mathbf{v}_p=0 \text{ on } \partial\Omega_p\},\\
X_p &= \{\psi\in H^1(\Omega_p) : \psi=0 \text{ on } \Gamma_{E}\},\\
T_p &= \{\omega\in H^1(\Omega_p) : \omega=0 \text{ on } \Gamma_{E}\}.
\end{align*}

From now on, the inner product in $L^2(D)$ ($D$ may be $\Omega_f$, $\Omega_p$ or $\Gamma$) is denoted by $(\cdot,\cdot)$, and the norm of $L^2(D)$ by $\|\cdot\|_D$. We omit the subscript $D$ when no ambiguity arises.

We derive the weak formulation of the coupled model \eqref{eq:fluid_momentum}-\eqref{eq:interface_heatflux} as follows: given
$\mathbf{f}_f \in L^2(\Omega_f)$,
$\mathbf{f}_p \in L^2(\Omega_p)$,
$g_p \in L^2(\Omega_p)$, find
$\mathbf{u}_f \in \mathbf{H}_f$,
$p_f \in X_f$,
$\theta_f \in T_f$,
$\etavec_p \in \mathbf{H}_p$,
$\xivec_p \in \mathbf{H}_p$,
$p_p \in X_p$,
$\theta_p \in T_p$
such that for any
$\mathbf{v}_f \in \mathbf{H}_f$,
$q_f \in X_f$,
$w_f \in T_f$,
$\bm{\varphi}_p \in \mathbf{H}_p$,
$q_p \in X_p$,
$w_p \in T_p$
and $t\in(0,T]$, there hold
\begin{align}
&\bigl(\partial_t\mathbf{u}_f,\mathbf{v}_f\bigr)_{\Omega_f} + 2\mu_f\bigl(\mathbb{D}(\mathbf{u}_f),\mathbb{D}(\mathbf{v}_f)\bigr)_{\Omega_f} + \langle p_p,\mathbf{v}_f\cdot\mathbf{n}_f\rangle_\Gamma + \frac{\alpha\sqrt{d}}{\sqrt{\operatorname{tr}(\mathbf{K})}}\bigl\langle P_{\boldsymbol{\tau}}(\mathbf{u}_f-\partial_t\etavec_p), P_{\boldsymbol{\tau}}(\mathbf{v}_f)\bigr\rangle_\Gamma\nonumber\\
&\quad  + \bigl((\mathbf{u}_f\cdot\nabla)\mathbf{u}_f,\mathbf{v}_f\bigr)_{\Omega_f}  - (\nabla\cdot\mathbf{v}_f,p_f)_{\Omega_f} = \frac{Gr}{Re^2}\bigl(\boldsymbol{\xi}\theta_f,\mathbf{v}_f\bigr)_{\Omega_f}, \label{eq:12}\\[2mm]
&(\nabla\cdot\mathbf{u}_f,q_f)_{\Omega_f} = 0, \label{eq:13}\\[2mm]
&\bigl(\partial_t\theta_f,w_f\bigr)_{\Omega_f} + a_f\bigl(\nabla\theta_f,\nabla w_f\bigr)_{\Omega_f} + \bigl((\mathbf{u}_f\cdot\nabla)\theta_f,w_f\bigr)_{\Omega_f}  - a_f\langle \nabla\theta_f\cdot\mathbf{n}_f,w_f\rangle_\Gamma = \bigl(\mathbf{u}_f\cdot\boldsymbol{\xi},w_f\bigr)_{\Omega_f}, \label{eq:14}\\[2mm]
&\bigl(\partial_t\xivec_p,\bm{\varphi}_p\bigr)_{\Omega_p} + 2\mu_p\bigl(\mathbb{D}(\etavec_p),\mathbb{D}(\bm{\varphi}_p)\bigr)_{\Omega_p} + \lambda_p(\nabla\cdot\etavec_p,\nabla\cdot\bm{\varphi}_p)_{\Omega_p}  - \alpha(\nabla\cdot\bm{\varphi}_p,p_p)_{\Omega_p} \nonumber\\
&\quad- \beta(\nabla\cdot\bm{\varphi}_p,\theta_p)_{\Omega_p}  - \langle p_p,\bm{\varphi}_p\cdot\mathbf{n}_f\rangle_\Gamma - \frac{\alpha\sqrt{d}}{\sqrt{\operatorname{tr}(\mathbf{K})}}\bigl\langle P_{\boldsymbol{\tau}}(\mathbf{u}_f-\partial_t\etavec_p), P_{\boldsymbol{\tau}}(\bm{\varphi}_p)\bigr\rangle_\Gamma = (\mathbf{f},\bm{\varphi}_p)_{\Omega_p}, \label{eq:15}\\[2mm]
&S_0\bigl(\partial_t p_p,q_p\bigr)_{\Omega_p} - b_0\bigl(\partial_t\theta_p,q_p\bigr)_{\Omega_p} + \alpha\bigl(\nabla\cdot\xivec_p,q_p\bigr)_{\Omega_p} + \mathbf{K}\bigl(\nabla p_p,\nabla q_p\bigr)_{\Omega_p}  \nonumber\\
&\quad- \langle (\mathbf{u}_f-\xivec_p)\cdot\mathbf{n}_f,q_p\rangle_\Gamma = (g,q_p)_{\Omega_p}, \label{eq:16}\\[2mm]
&a_0\bigl(\partial_t\theta_p,w_p\bigr)_{\Omega_p} - b_0\bigl(\partial_t p_p,w_p\bigr)_{\Omega_p} + \beta\bigl(\nabla\cdot\xivec_p,w_p\bigr)_{\Omega_p} + \bm{a}_p\bigl(\nabla\theta_p,\nabla w_p\bigr)_{\Omega_p}  \nonumber\\
&\quad- \bm{a}_p\langle \nabla\theta_p\cdot\mathbf{n}_p,w_p\rangle_\Gamma = (\phi,w_p)_{\Omega_p},\label{eq:17}\\
&\mathbf{u}_f(0)=\mathbf{0},\,\boldsymbol{\eta}_p(0)=\mathbf{0},\,\boldsymbol{\xi}_p(0)=\mathbf{0},\, p_f(0)=0,\, p_p(0)=0,\,\theta_f(0)=0,\, \theta_p(0)=0
. \label{eq:18}
\end{align}

\begin{remark}\cite{chen2022multiphysics}
The following assumptions are imposed on the problem data.
\begin{itemize}
    \item The tensor $\mathbf{K} : \mathbb{R}^d \to \mathbb{R}^{d \times d}$ is symmetric and uniformly positive definite, i.e., there exist positive constants $K_{min}$ and $K_{max}$ such that
    \begin{align*}
    K_{min} |\zeta|^2 \le \zeta^\mathrm{T} \mathbf{K}(x) \zeta \le K_{max} |\zeta|^2, \qquad \forall \zeta \in \mathbb{R}^d \setminus \{0\}.
    \end{align*}

    \item The tensor $\bm{a}_p : \mathbb{R}^d \to \mathbb{R}^{d \times d}$ is symmetric and uniformly positive definite, i.e., there exist positive constants ${a_p}_{min}$ and ${a_p}_{max}$ such that
    \begin{align*}
    {a_p}_{min} |\zeta|^2 \le \zeta^\mathrm{T} \bm{a}_p(x) \zeta \le {a_p}_{max} |\zeta|^2, \qquad \forall \zeta \in \mathbb{R}^d \setminus \{0\}.
   \end{align*}

    \item The coefficients $S_0$, $a_0$ and $b_0$ are nonnegative constants satisfying $S_0 - b_0 > 0$ and $a_0 - b_0 > 0$.
\end{itemize}
\end{remark}

We next recall two classical results that will play a crucial role in the analysis.
\begin{lemma}\cite{cesmelioglu2017analysis}[Gronwall's inequality in integral form]\label{thm:gronwall}
Let $\gamma$ be a continuous nonnegative function and let $B, C \ge 0$ be constants such that
\begin{align*}
\gamma(t) \le B + C \int_0^t \gamma(s)\,ds \qquad \text{a.e. } t \in (0,T).
\end{align*}
Then $\gamma(t) \le B e^{Ct}$ for all $t \in [0,T]$.
\end{lemma}

In the subsequent passage to the limit, we will encounter the convergence issue of the nonlinear convective terms. This is achieved by invoking the following compactness lemma, which guarantees strong convergence in appropriate spaces.
\begin{lemma}\cite{simon1986compact}\label{lem:compactness}
Let $X$, $B$ and $Y$ be Banach spaces such that $X \subset B \subset Y$ where the embedding of $X$ into $B$ is compact. Let $F$ be a bounded set in $L^p(0,T;X)$ where $1 \le p < \infty$ and let the set $\bigl\{ \frac{\partial f}{\partial t} \bigr\}_{f \in F}$ be bounded in $L^1(0,T;Y)$. Then $F$ is relatively compact in $L^p(0,T;B)$.
\end{lemma}

\section{Existence analysis}\label{sec:3}
In this section, we prove the existence of a weak solution to problem \eqref{eq:12}--\eqref{eq:18}. The proof follows a classical Galerkin method and consists of several steps. First, we introduce a divergence-free weak formulation where the fluid pressure is eliminated. A sequence of Galerkin approximations is then constructed, uniform a priori estimates are derived, and the limit is passed to the infinite-dimensional system.
\subsection{Galerkin approximations}

To construct approximate solutions via the Galerkin method, it is convenient to first eliminate the fluid pressure $p_f$ from the weak formulation. This is achieved by restricting the test functions for the fluid momentum equation to the space of divergence-free functions, which we define as
\begin{align*}
\mathbf{V}_f := \{\mathbf{v}\in \mathbf{H}_f : \nabla\cdot\mathbf{v}=0\}.
\end{align*}

Consequently, the pressure term $(\nabla\cdot\mathbf{v}_f,p_f)$
vanishes, and the pressure
$p_f$
  no longer appears in the formulation. The resulting system, employed for the Galerkin approximation, is given by:
given
$\mathbf{f}_f \in L^2(\Omega_f)$,
$\mathbf{f}_p \in L^2(\Omega_p)$,
$g_p \in L^2(\Omega_p)$, find
$\mathbf{u}_f \in \mathbf{V}_f$,
$\theta_f \in T_f$,
$\etavec_p \in \mathbf{H}_p$,
$\xivec_p \in \mathbf{H}_p$,
$p_p \in X_p$,
$\theta_p \in T_p$
such that for any
$\mathbf{v}_f \in \mathbf{V}_f$,
$w_f \in T_f$,
$\bm{\varphi}_p \in \mathbf{H}_p$,
$q_p \in X_p$,
$w_p \in T_p$
and $t\in(0,T]$, there hold
\begin{align}
&\bigl(\partial_t\mathbf{u}_f,\mathbf{v}_f\bigr)_{\Omega_f} + 2\mu_f\bigl(\mathbb{D}(\mathbf{u}_f),\mathbb{D}(\mathbf{v}_f)\bigr)_{\Omega_f} + \langle p_p,\mathbf{v}_f\cdot\mathbf{n}_f\rangle_\Gamma + \frac{\alpha\sqrt{d}}{\sqrt{\operatorname{tr}(\mathbf{K})}}\bigl\langle P_{\boldsymbol{\tau}}(\mathbf{u}_f-\partial_t\etavec_p), P_{\boldsymbol{\tau}}(\mathbf{v}_f)\bigr\rangle_\Gamma\nonumber\\
&\quad  + \bigl((\mathbf{u}_f\cdot\nabla)\mathbf{u}_f,\mathbf{v}_f\bigr)_{\Omega_f}   = \frac{Gr}{Re^2}\bigl(\boldsymbol{\xi}\theta_f,\mathbf{v}_f\bigr)_{\Omega_f}, \label{eq:19}\\[2mm]
&(\nabla\cdot\mathbf{u}_f,q_f)_{\Omega_f} = 0, \label{eq:20}\\[2mm]
&\bigl(\partial_t\theta_f,w_f\bigr)_{\Omega_f} + a_f\bigl(\nabla\theta_f,\nabla w_f\bigr)_{\Omega_f} + \bigl((\mathbf{u}_f\cdot\nabla)\theta_f,w_f\bigr)_{\Omega_f}  - a_f\langle \nabla\theta_f\cdot\mathbf{n}_f,w_f\rangle_\Gamma = \bigl(\mathbf{u}_f\cdot\boldsymbol{\xi},w_f\bigr)_{\Omega_f}, \label{eq:21}\\[2mm]
&\bigl(\partial_t\xivec_p,\bm{\varphi}_p\bigr)_{\Omega_p} + 2\mu_p\bigl(\mathbb{D}(\etavec_p),\mathbb{D}(\bm{\varphi}_p)\bigr)_{\Omega_p} + \lambda_p(\nabla\cdot\etavec_p,\nabla\cdot\bm{\varphi}_p)_{\Omega_p}  - \alpha(\nabla\cdot\bm{\varphi}_p,p_p)_{\Omega_p} \nonumber\\
&\quad- \beta(\nabla\cdot\bm{\varphi}_p,\theta_p)_{\Omega_p}  - \langle p_p,\bm{\varphi}_p\cdot\mathbf{n}_f\rangle_\Gamma - \frac{\alpha\sqrt{d}}{\sqrt{\operatorname{tr}(\mathbf{K})}}\bigl\langle P_{\boldsymbol{\tau}}(\mathbf{u}_f-\partial_t\etavec_p), P_{\boldsymbol{\tau}}(\bm{\varphi}_p)\bigr\rangle_\Gamma = (\mathbf{f},\bm{\varphi}_p)_{\Omega_p}, \label{eq:22}\\[2mm]
&S_0\bigl(\partial_t p_p,q_p\bigr)_{\Omega_p} - b_0\bigl(\partial_t\theta_p,q_p\bigr)_{\Omega_p} + \alpha\bigl(\nabla\cdot\xivec_p,q_p\bigr)_{\Omega_p} + \mathbf{K}\bigl(\nabla p_p,\nabla q_p\bigr)_{\Omega_p}  \nonumber\\
&\quad- \langle (\mathbf{u}_f-\xivec_p)\cdot\mathbf{n}_f,q_p\rangle_\Gamma = (g,q_p)_{\Omega_p}, \label{eq:23}\\[2mm]
&a_0\bigl(\partial_t\theta_p,w_p\bigr)_{\Omega_p} - b_0\bigl(\partial_t p_p,w_p\bigr)_{\Omega_p} + \beta\bigl(\nabla\cdot\xivec_p,w_p\bigr)_{\Omega_p} + \bm{a}_p\bigl(\nabla\theta_p,\nabla w_p\bigr)_{\Omega_p}  \nonumber\\
&\quad- \bm{a}_p\langle \nabla\theta_p\cdot\mathbf{n}_p,w_p\rangle_\Gamma = (\phi,w_p)_{\Omega_p},\label{eq:24}
\end{align}
and
\begin{align}
\mathbf{u}_f(0)=\mathbf{0},\, \theta_f(0)=0\,\text{a.e. in }\Omega_f \quad\boldsymbol{\eta}_p(0)=\mathbf{0},\,\boldsymbol{\xi}_p(0)=\mathbf{0},\, p_p(0)=0,\, \theta_p(0)=0\,\text{a.e. in }\Omega_p.\label{eq:25}
\end{align}
With $p_f$ eliminated, it is clear that every solution to \eqref{eq:12}--\eqref{eq:18} also satisfies \eqref{eq:19}--\eqref{eq:25}.

To construct approximate solutions, we apply the Galerkin method based on finite-dimensional subspaces spanned by basis functions of the separable spaces introduced above.

For each of the separable spaces \(\mathbf{V}_f, T_f, \mathbf{H}_p, X_p, T_p\), choose a smooth dense basis, denoted respectively by
\[
\{\mathbf{v}_{f,i}\},\ \{w_{f,i}\},\ \{\boldsymbol{\varphi}_{p,i}\},\ \{q_{p,i}\},\ \{w_{p,i}\}.
\]
Given \(m\in\mathbb{N}\), denote by \(\mathbf{V}_f^m, T_f^m, \mathbf{H}_p^m, X_p^m, T_p^m\) the spans of the first \(m\) elements of the corresponding bases, and use the following Galerkin approximations for the unknowns \(\mathbf{u}_{f}\), \(\theta_{f}\), \(\boldsymbol{\eta}_{p}\), \(p_{p}\) and \(\theta_{p}\):
\begin{equation}
\begin{aligned}
\mathbf{u}_{f,m}(\mathbf{x},t) &= \sum_{j=1}^m a_j(t)\,\mathbf{v}_{f,j}(\mathbf{x}), &
\theta_{f,m}(\mathbf{x},t) &= \sum_{j=1}^m b_j(t)\,w_{f,j}(\mathbf{x}),
&\boldsymbol{\eta}_{p,m}(\mathbf{x},t) &= \sum_{j=1}^m c_j(t)\,\boldsymbol{\varphi}_{p,j}(\mathbf{x}), \\
p_{p,m}(\mathbf{x},t) &= \sum_{j=1}^m d_j(t)\,q_{p,j}(\mathbf{x}), &
\theta_{p,m}(\mathbf{x},t) &= \sum_{j=1}^m e_j(t)\,w_{p,j}(\mathbf{x}).\label{eq:26}
\end{aligned}
\end{equation}
We note that \(\boldsymbol{\xi}_{p,m}:=\partial_t\boldsymbol{\eta}_{p,m}\), so that
\begin{align}
\boldsymbol{\xi}_{p,m}(\mathbf{x},t) = \sum_{j=1}^m \partial_t {c}_j(t)\,\boldsymbol{\varphi}_{p,j}(\mathbf{x}),\label{eq:27}
\end{align}
where \(\partial_t {c}_j(t)\) denotes the derivative of \(c_j(t)\) with respect to \(t\).

Insert \eqref{eq:26} and \eqref{eq:27} into the weak formulation \eqref{eq:19}--\eqref{eq:25}. For each $i=1,\dots,m$ we obtain the following five equations:
\begin{align}\label{eq:28}
&\bigl(\partial_t\mathbf{u}_{f,m},\mathbf{v}_{f}\bigr)_{\Omega_f}
+2\mu_f\bigl(\mathbb{D}(\mathbf{u}_{f,m}),\mathbb{D}(\mathbf{v}_{f})\bigr)_{\Omega_f}
+\ip{p_{p,m}}{\mathbf{v}_{f}\!\cdot\!\mathbf{n}_f}_\Gamma+\bigl((\mathbf{u}_{f,m}\!\cdot\!\nabla)\mathbf{u}_{f,m},\mathbf{v}_{f}\bigr)_{\Omega_f}
\nonumber\\
&\quad+\frac{\alpha\sqrt{d}}{\sqrt{\operatorname{tr}(\mathbf{K})}}\ip{P_\tau(\mathbf{u}_{f,m}-\boldsymbol{\xi}_{p,m})}{P_\tau(\mathbf{v}_{f})}_\Gamma
= \frac{G_r}{Re^2}\bigl(\bm{\xi}\theta_{f,m},\mathbf{v}_{f}\bigr)_{\Omega_f},\\ \label{eq:29}
&\bigl(\partial_t\theta_{f,m},w_{f}\bigr)_{\Omega_f}
+a_f\bigl(\nabla\theta_{f,m},\nabla w_{f}\bigr)_{\Omega_f}
+\bigl((\mathbf{u}_{f,m}\!\cdot\!\nabla)\theta_{f,m},w_{f}\bigr)_{\Omega_f}
\nonumber\\
&\quad-a_f\ip{\nabla\theta_{f,m}\!\cdot\!\mathbf{n}_f}{w_{f}}_\Gamma
= \bigl(\mathbf{u}_{f,m}\!\cdot\!\bm{\xi},w_{f}\bigr)_{\Omega_f},\\\label{eq:30}
&\bigl(\partial_t\bm{\xi}_{p,m},\bm{\varphi}_{p}\bigr)_{\Omega_p}
+2\mu_p\bigl(\mathbb{D}(\bm{\eta}_{p,m}),\mathbb{D}(\bm{\varphi}_{p})\bigr)_{\Omega_p}
+\lambda_p\bigl(\nabla\!\cdot\!\boldsymbol{\eta}_{p,m},\nabla\!\cdot\!\bm{\varphi}_{p}\bigr)_{\Omega_p}
-\alpha\bigl(\nabla\!\cdot\!\bm{\varphi}_{p},p_{p,m}\bigr)_{\Omega_p}
\nonumber\\
&\quad-\beta\bigl(\nabla\!\cdot\!\bm{\varphi}_{p},\theta_{p,m}\bigr)_{\Omega_p}
-\ip{p_{p,m}}{\bm{\varphi}_{p}\!\cdot\!\mathbf{n}_f}_\Gamma
-\frac{\alpha\sqrt{d}}{\sqrt{\operatorname{tr}(\mathbf{K})}}\ip{P_\tau(\mathbf{u}_{f,m}-\boldsymbol{\xi}_{p,m})}{P_\tau(\bm{\varphi}_{p})}_\Gamma
= \bigl(\mathbf{f},\bm{\varphi}_{p}\bigr)_{\Omega_p},\\\label{eq:31}
&S_0\bigl(\partial_t p_{p,m},q_{p}\bigr)_{\Omega_p}
-b_0\bigl(\partial_t\theta_{p,m},q_{p}\bigr)_{\Omega_p}
+\alpha\bigl(\nabla\!\cdot\!\boldsymbol{\xi}_{p,m},q_{p}\bigr)_{\Omega_p}
+\mathbf{K}\bigl(\nabla p_{p,m},\nabla q_{p}\bigr)_{\Omega_p}\nonumber\\
&\quad -\ip{(\mathbf{u}_{f,m}-\boldsymbol{\xi}_{p,m})\!\cdot\!\mathbf{n}_f}{q_{p}}_\Gamma
= \bigl(g,q_{p}\bigr)_{\Omega_p}.\\
&a_0\bigl(\partial_t\theta_{p,m},w_{p}\bigr)_{\Omega_p}
-b_0\bigl(\partial_t p_{p,m},w_{p}\bigr)_{\Omega_p}
+\beta\bigl(\nabla\!\cdot\!\boldsymbol{\xi}_{p,m},w_{p}\bigr)_{\Omega_p}
+\bm{a}_p\bigl(\nabla\theta_{p,m},\nabla w_{p}\bigr)_{\Omega_p}
\nonumber\\
&\quad-\bm{a}_p\ip{\nabla\theta_{p,m}\!\cdot\!\mathbf{n}_p}{w_{p}}_\Gamma = \bigl(\phi,w_{p}\bigr)_{\Omega_p}.\label{eq:32}
\end{align}
and
\begin{align}
\mathbf{u}_{f,m}(0)=\mathbf{0},\, \theta_{f,m}(0)=0,\, \boldsymbol{\eta}_{p,m}(0)=\mathbf{0},\,\boldsymbol{\xi}_{p,m}(0)=\mathbf{0},\, p_{p,m}(0)=0,\, \theta_{p,m}(0)=0.\label{eq:33}
\end{align}

\label{subsec:galerkin}
\begin{theorem}
For each positive integer $m$, equations \eqref{eq:28}-\eqref{eq:33} have a unique maximal solution
$\mathbf{u}_{f,m} \in C^1(0, T_m; \mathbf{V}^m_f)$, $\theta_{f,m} \in C^1(0, T_m; T^m_f)$,
$\boldsymbol{\eta}_{p,m} \in C^2(0, T_m; \mathbf{H}^m_p)$,
$\boldsymbol{\xi}_{p,m} \in C^1(0, T_m; \mathbf{H}^m_p)$,
$p_{p,m} \in C^1(0, T_m; X^m_p)$,
$\theta_{p,m} \in C^1(0, T_m; T^m_p)$
for some time $T_m$ where $0 < T_m \leq T$.
\end{theorem}
\begin{proof}\
Using the Galerkin expansions given in \eqref{eq:26} and \eqref{eq:27}, and choosing the basis functions as the test functions. For each $i=1,\dots,m$ we obtain:
\begin{align}\label{eq:34}
&\sum_{j=1}^m \partial_t a_j ( \mathbf{v}_{f,j},\mathbf{v}_{f,i})_{\Omega_f}
+ 2\mu_f \sum_{j=1}^m a_j (\mathbb{D}(\mathbf{v}_{f,j}),\mathbb{D}(\mathbf{v}_{f,i}))_{\Omega_f}
+ \sum_{j=1}^m d_j \langle q_{p,j},\mathbf{v}_{f,i}\cdot\mathbf{n}_f\rangle_{\Gamma} \nonumber\\
&+ \frac{\alpha\sqrt{d}}{\sqrt{\operatorname{tr}(\mathbf{K})}}
\sum_{j=1}^m \bigl[ a_j \langle P_\tau(\mathbf{v}_{f,j}),P_\tau(\mathbf{v}_{f,i})\rangle_{\Gamma}
- \partial_t c_j \langle P_\tau(\boldsymbol{\varphi}_{p,j}),P_\tau(\mathbf{v}_{f,i})\rangle_{\Gamma} \bigr] \nonumber\\
&+ \sum_{j,k=1}^m a_j a_k \bigl((\mathbf{v}_{f,j}\cdot\nabla)\mathbf{v}_{f,k},\mathbf{v}_{f,i}\bigr)_{\Omega_f}
= \frac{Gr}{Re^2} \sum_{j=1}^m b_j (\boldsymbol{\xi} w_{f,j},\mathbf{v}_{f,i})_{\Omega_f}, \\
&\sum_{j=1}^m \partial_t b_j ( w_{f,j},w_{f,i})_{\Omega_f}
+ a_f \sum_{j=1}^m b_j (\nabla w_{f,j},\nabla w_{f,i})_{\Omega_f}
+ \sum_{j,k=1}^m a_j b_k \bigl((\mathbf{v}_{f,j}\cdot\nabla)w_{f,k},w_{f,i}\bigr)_{\Omega_f} \nonumber\\
&- a_f \sum_{j=1}^m b_j \langle \nabla w_{f,j}\cdot\mathbf{n}_f,w_{f,i}\rangle_{\Gamma}
= \sum_{j=1}^m a_j (\mathbf{v}_{f,j}\cdot\boldsymbol{\xi},w_{f,i})_{\Omega_f}, \\
&\sum_{j=1}^m \partial_t^2 c_j ( \boldsymbol{\varphi}_{p,j},\boldsymbol{\varphi}_{p,i})_{\Omega_p}
+ 2\mu_p \sum_{j=1}^m c_j (\mathbb{D}(\boldsymbol{\varphi}_{p,j}),\mathbb{D}(\boldsymbol{\varphi}_{p,i}))_{\Omega_p}
+ \lambda_p \sum_{j=1}^m c_j (\nabla\!\cdot\!\boldsymbol{\varphi}_{p,j},\nabla\!\cdot\!\boldsymbol{\varphi}_{p,i})_{\Omega_p} \nonumber\\
&- \alpha \sum_{j=1}^m d_j (\nabla\!\cdot\!\boldsymbol{\varphi}_{p,i},q_{p,j})_{\Omega_p}
- \beta \sum_{j=1}^m e_j (\nabla\!\cdot\!\boldsymbol{\varphi}_{p,i},w_{p,j})_{\Omega_p}
- \sum_{j=1}^m d_j \langle q_{p,j},\boldsymbol{\varphi}_{p,i}\cdot\mathbf{n}_f\rangle_{\Gamma} \nonumber\\
&- \frac{\alpha\sqrt{d}}{\sqrt{\operatorname{tr}(\mathbf{K})}}
\sum_{j=1}^m \bigl[ a_j \langle P_\tau(\mathbf{v}_{f,j}),P_\tau(\boldsymbol{\varphi}_{p,i})\rangle_{\Gamma}
- \partial_t c_j \langle P_\tau(\boldsymbol{\varphi}_{p,j}),P_\tau(\boldsymbol{\varphi}_{p,i})\rangle_{\Gamma} \bigr]
= (\mathbf{f},\boldsymbol{\varphi}_{p,i})_{\Omega_p}, \\
&S_0 \sum_{j=1}^m \partial_t d_j ( q_{p,j},q_{p,i})_{\Omega_p}
- b_0 \sum_{j=1}^m \partial_t e_j ( w_{p,j},q_{p,i})_{\Omega_p}
+ \alpha \sum_{j=1}^m \partial_t c_j (\nabla\!\cdot\!\boldsymbol{\varphi}_{p,j},q_{p,i})_{\Omega_p} \nonumber\\
&+ \sum_{j=1}^m d_j (\mathbf{K}\nabla q_{p,j},\nabla q_{p,i})_{\Omega_p}
- \sum_{j=1}^m \bigl[ a_j \langle \mathbf{v}_{f,j}\cdot\mathbf{n}_f,q_{p,i}\rangle_{\Gamma}
- \partial_t c_j \langle \boldsymbol{\varphi}_{p,j}\cdot\mathbf{n}_f,q_{p,i}\rangle_{\Gamma} \bigr]
= (g,q_{p,i})_{\Omega_p}, \\
&a_0 \sum_{j=1}^m \partial_t e_j ( w_{p,j},w_{p,i})_{\Omega_p}
- b_0 \sum_{j=1}^m \partial_t d_j ( q_{p,j},w_{p,i})_{\Omega_p}
+ \beta \sum_{j=1}^m \partial_t c_j (\nabla\!\cdot\!\boldsymbol{\varphi}_{p,j},w_{p,i})_{\Omega_p} \nonumber\\
&+ \bm{a}_p \sum_{j=1}^m e_j (\nabla w_{p,j},\nabla w_{p,i})_{\Omega_p}
+ \bm{a}_p \sum_{j=1}^m e_j \langle \nabla w_{p,j}\cdot\mathbf{n}_f,w_{p,i}\rangle_{\Gamma}
= (\phi,w_{p,i})_{\Omega_p}.\label{eq:38}
\end{align}
Define the vectors of coefficients
\begin{align*}
&\mathbf{a}(t)=(a_1,\dots,a_m)^\mathrm{T},\quad \mathbf{b}(t)=(b_1,\dots,b_m)^\mathrm{T},\quad \mathbf{c}(t)=(c_1,\dots,c_m)^\mathrm{T},\\
&\mathbf{d}(t)=(d_1,\dots,d_m)^\mathrm{T},\quad \mathbf{e}(t)=(e_1,\dots,e_m)^\mathrm{T},
\end{align*}
and introduce the extended state vector
\begin{align*}
\mathbf{W}(t)=\bigl(\mathbf{a}(t),\mathbf{b}(t),\partial_t {\mathbf{c}}(t),\mathbf{c}(t),\mathbf{d}(t),\mathbf{e}(t)\bigr)^\mathrm{T}\in\mathbb{R}^{6m}.
\end{align*}
By employing the Galerkin expansions from \eqref{eq:26} and \eqref{eq:27}, problem \eqref{eq:34}--\eqref{eq:38} can be rewritten in matrix form. This leads to a standard finite-dimensional argument, effectively formulating the problem as a square first‑order system of ordinary differential equations with an initial condition. For the integrals appearing on the left‑hand side, with
$1 \le i,\; j \le m$, we define:
\begin{gather*}
(M_f)_{ij}=(\mathbf{v}_{f,j},\mathbf{v}_{f,i})_{\Omega_f},\quad
(M_{Tf})_{ij}=(w_{f,j},w_{f,i})_{\Omega_f},\quad
(M_s)_{ij}=(\boldsymbol{\varphi}_{p,j},\boldsymbol{\varphi}_{p,i})_{\Omega_p},\\[2mm]
(M_p)_{ij}=S_0(q_{p,j},q_{p,i})_{\Omega_p},\quad
(M_{Tp})_{ij}=a_0(w_{p,j},w_{p,i})_{\Omega_p},\quad
(N_{pT})_{ij}=b_0(w_{p,j},q_{p,i})_{\Omega_p},\\[2mm]
(A_f)_{ij}=2\mu_f(\mathbb{D}(\mathbf{v}_{f,j}),\mathbb{D}(\mathbf{v}_{f,i}))_{\Omega_f},\quad
(A_{Tf})_{ij}=a_f(\nabla w_{f,j},\nabla w_{f,i})_{\Omega_f},\\[2mm]
(A_s)_{ij}=2\mu_p(\mathbb{D}(\boldsymbol{\varphi}_{p,j}),\mathbb{D}(\boldsymbol{\varphi}_{p,i}))_{\Omega_p}
+\lambda_p(\nabla\!\cdot\!\boldsymbol{\varphi}_{p,j},\nabla\!\cdot\!\boldsymbol{\varphi}_{p,i})_{\Omega_p},\\[2mm]
(A_p)_{ij}=(\mathbf{K}\nabla q_{p,j},\nabla q_{p,i})_{\Omega_p},\quad
(A_{Tp})_{ij}=\bm{a}_p(\nabla w_{p,j},\nabla w_{p,i})_{\Omega_p},\\[2mm]
(B_{ff})_{ij}=\frac{\alpha\sqrt{d}}{\sqrt{\operatorname{tr}(\mathbf{K})}}\langle P_\tau(\mathbf{v}_{f,j}),P_\tau(\mathbf{v}_{f,i})\rangle_\Gamma,\quad
(B_{fs})_{ij}=-\frac{\alpha\sqrt{d}}{\sqrt{\operatorname{tr}(\mathbf{K})}}\langle P_\tau(\boldsymbol{\varphi}_{p,j}),P_\tau(\mathbf{v}_{f,i})\rangle_\Gamma,\\[2mm]
(B_{ss})_{ij}=\frac{\alpha\sqrt{d}}{\sqrt{\operatorname{tr}(\mathbf{K})}}\langle P_\tau(\boldsymbol{\varphi}_{p,j}),P_\tau(\boldsymbol{\varphi}_{p,i})\rangle_\Gamma,\quad
(C_{pf})_{ij}=\langle q_{p,j},\mathbf{v}_{f,i}\cdot\mathbf{n}_f\rangle_\Gamma,\\[2mm]
(C_{ps})_{ij}=\alpha(\nabla\!\cdot\!\boldsymbol{\varphi}_{p,i},q_{p,j})_{\Omega_p}+\langle q_{p,j},\boldsymbol{\varphi}_{p,i}\cdot\mathbf{n}_f\rangle_\Gamma,\quad
(C_{Ts})_{ij}=\beta(\nabla\!\cdot\!\boldsymbol{\varphi}_{p,i},w_{p,j})_{\Omega_p},\\[2mm]
(L_{Tf})_{ij}=\langle \nabla a_fw_{f,j}\cdot\mathbf{n}_f,w_{f,i}\rangle_\Gamma,\quad
(D_s)_{ij}=\langle \boldsymbol{\varphi}_{p,j}\cdot\mathbf{n}_f,q_{p,i}\rangle_\Gamma,\quad
(Q_s)_{ij}=\alpha(\nabla\!\cdot\!\boldsymbol{\varphi}_{p,j},q_{p,i})_{\Omega_p},\\[2mm]
(Q_{Ts})_{ij}=\beta(\nabla\!\cdot\!\boldsymbol{\varphi}_{p,j},w_{p,i})_{\Omega_p},\quad
(R_{Tp})_{ij}=\bm{a}_p\langle \nabla w_{p,j}\cdot\mathbf{n}_f,w_{p,i}\rangle_\Gamma,\quad
(C_{fT})_{ij}=\frac{Gr}{Re^2}(\boldsymbol{\xi} w_{f,j},\mathbf{v}_{f,i})_{\Omega_f},\\[2mm]
(C_{Tf})_{ij}=(\mathbf{v}_{f,j}\cdot\boldsymbol{\xi},w_{f,i})_{\Omega_f}.
\end{gather*}
Define the nonlinear vectors
\begin{align*}
(\mathbf{N}_{\text{conv}}(\mathbf{a}))_i &= \sum_{j,k=1}^m a_j a_k \bigl((\mathbf{v}_{f,j}\cdot\nabla)\mathbf{v}_{f,k},\mathbf{v}_{f,i}\bigr)_{\Omega_f},\\
(\mathbf{N}_{\text{therm}}(\mathbf{a},\mathbf{b}))_i &= \sum_{j,k=1}^m a_j b_k \bigl((\mathbf{v}_{f,j}\cdot\nabla)w_{f,k},w_{f,i}\bigr)_{\Omega_f}.
\end{align*}
The right-hand side vectors are
\begin{align*}
(\mathbf{F}_s)_i = (\mathbf{f},\boldsymbol{\varphi}_{p,i})_{\Omega_p},\quad
(\mathbf{G}_p)_i = (g,q_{p,i})_{\Omega_p},\quad
(\boldsymbol{\Phi})_i = (\phi,w_{p,i})_{\Omega_p}.
\end{align*}
Using these definitions, the system \eqref{eq:28}--\eqref{eq:33} can be written compactly as
\begin{align*}
\mathcal{M}\,\frac{d}{dt}\mathbf{W} + \mathcal{L}\,\mathbf{W} + \mathcal{N}(\mathbf{W}) = \mathbf{G}(t),
\end{align*}
where
\begin{align*}
\mathcal{M}=
\begin{pmatrix}
M_f & 0 & 0 & 0 & 0 & 0\\
0 & M_{Tf} & 0 & 0 & 0 & 0\\
0 & 0 & M_s & 0 & 0 & 0\\
0 & 0 & 0 & I_m & 0 & 0\\
0 & 0 & 0 & 0 & M_p & -N_{pT}\\
0 & 0 & 0 & 0 & -N_{pT}^\mathrm{T} & M_{Tp}
\end{pmatrix},
\end{align*}

\noindent \textbf{Linear stiffness matrix $\mathcal{L}$}\\
\begin{align*}
\mathcal{L}=
\begin{pmatrix}
A_f + B_{ff}-C_{fT}& 0 & -B_{fs} & 0 & C_{pf} & 0\\
-C_{Tf} & A_{Tf} - L_{Tf} & 0 & 0 & 0 & 0\\
-B_{fs}^\mathrm{T} & 0 & B_{ss} & 0 & 0 & 0\\
0 & 0 & 0 & A_s & -C_{ps}^\mathrm{T} & -C_{Ts}^\mathrm{T}\\
-D_s^T & 0 & Q_s + D_s & 0 & A_p & 0\\
0 & 0 & Q_{Ts} & 0 & 0 & A_{Tp} + R_{Tp}
\end{pmatrix},
\end{align*}

\noindent \textbf{Nonlinear term $\mathcal{N}(\mathbf{W})$}\\
\begin{align*}
\mathcal{N}(\mathbf{W}) = \bigl( \mathbf{N}_{\text{conv}}(\mathbf{a}),\; \mathbf{N}_{\text{therm}}(\mathbf{a},\mathbf{b}),\; \mathbf{0},\; \mathbf{0},\; \mathbf{0},\; \mathbf{0} \bigr)^\mathrm{T},
\end{align*}
\noindent \textbf{Right-hand side $\mathbf{G}(t,\mathbf{W})$ }\\
\begin{align*}
\mathbf{G}(t,\mathbf{W}) =
\begin{pmatrix}
\mathbf{0}, & \mathbf{0}, & \mathbf{F}_s(t), & \mathbf{0}, & \mathbf{G}_p(t), & \boldsymbol{\Phi}(t)
\end{pmatrix}^\mathrm{T}.
\end{align*}

All diagonal blocks of $\mathcal{M}$ are Gram matrices of linearly independent sets, hence symmetric positive definite. The $2\times2$ block
\[
\begin{pmatrix} M_p & -N_{pT} \\ -N_{pT}^\mathrm{T} & M_{Tp} \end{pmatrix}
\]
is positive definite because the physical parameters satisfy $S_0a_0 > b_0^2$ (thermo-poroelastic stability condition). Therefore $\mathcal{M}$ is invertible.

The  $ \mathcal{L}\mathbf{W} + \mathcal{N}(\mathbf{W})$ is a polynomial (linear plus quadratic), thus locally Lipschitz. The right-hand side $\mathbf{G}(t,\mathbf{W})$ is linear in $\mathbf{W}$ and continuous in $t$ by the regularity of the data. Consequently, the system can be rewritten as an explicit first-order ODE
\[
\frac{d}{dt}\mathbf{W} = \mathcal{M}^{-1}\bigl( \mathbf{G}(t,\mathbf{W}) - \mathcal{L}\mathbf{W} - \mathcal{N}(\mathbf{W}) \bigr) \equiv \mathbf{H}(t,\mathbf{W}),
\]
where $\mathbf{H}$ is continuous in $t$ and locally Lipschitz in $\mathbf{W}$.
With the initial condition $\mathbf{W}(0)=\mathbf{0}$, the Picard--Lindel\"of theorem
guarantees a unique local solution on some interval $[0,h]$.
By the standard continuation argument, this local solution can be extended
to a unique maximal solution $\mathbf{W}(t)$ on $[0,T_m)$ with $0<T_m\le T$,
and $\mathbf{W}\in C^1([0,T_m))$.
\end{proof}
\begin{remark}[Picard--Lindel\"of theorem]
Consider the initial value problem \(y'(t) = f(t, y(t))\), \(y(t_0) = y_0\).
If \(f\) is continuous and Lipschitz continuous with respect to \(y\) on a
closed rectangle \(R\) containing \((t_0, y_0)\) in its interior, then the
problem admits a unique local solution on some interval \([t_0 - h, t_0 + h]\).
\end{remark}

The remainder of the analysis concerns the derivation of a priori estimates uniform in \(m\), which show that \(T_m=T\), so that the Galerkin solution exists on the whole interval \([0,T]\). Subsequently, the limit \(m\to\infty\) is taken.

\subsection{A priori estimates for the Galerkin solution}
We present the inequalities and results that will be employed throughout the paper. Further details are available in \cite{ding1996proof,ern2004theory}.
\begin{itemize}
    \item \textbf{Poincaré inequality:} There exist constants $C_{P,f}$ and $C_{P,p}$ such that for every
          $\mathbf{v}\in\mathbf{H}_f$, $\theta\in T_f$, $q\in X_p$, $\omega\in T_p$,
          \begin{align}
              \norm{\mathbf{v}}_{L^2(\Omf)} &\le C_{P,f} \norm{\nabla\mathbf{v}}_{L^2(\Omf)}, &
              \norm{\theta}_{L^2(\Omf)} &\le C_{P,f} \norm{\nabla\theta}_{L^2(\Omf)}, \\
              \norm{q}_{L^2(\Omp)} &\le C_{P,p} \norm{\nabla q}_{L^2(\Omp)}, &
              \norm{\omega}_{L^2(\Omp)} &\le C_{P,p} \norm{\nabla\omega}_{L^2(\Omp)}.
          \end{align}

    \item \textbf{Sobolev embedding $H^1\hookrightarrow L^4$:} There exists a constant $C_{S,f}$ such that for every
          $\mathbf{v}\in\mathbf{H}_f$, $\theta\in T_f$,
          \begin{align}
              \norm{\mathbf{v}}_{L^4(\Omf)} &\le C_{S,f} \norm{\mathbf{v}}_{H^1(\Omf)}, \\
              \norm{\theta}_{L^4(\Omf)} &\le C_{S,f} \norm{\theta}_{H^1(\Omf)}.
          \end{align}

    \item \textbf{Korn's inequality:} There exists a constant $C_K>0$ such that for every $\mathbf{v}\in\mathbf{H}_f$,
          \begin{equation}
              \norm{\nabla\mathbf{v}}_{L^2} \le C_K \norm{\D(\mathbf{v})}_{L^2}.
          \end{equation}

    \item \textbf{Trace inequalities:} There exist constants $T_1, T_3, T_5$ such that
          \begin{align}
              \norm{\mathbf{v}}_{L^2(\Gamma)} &\le T_1 \norm{\mathbf{v}}_{H^1(\Omega_f)} \quad \forall \mathbf{v}\in\mathbf{H}_f, \\
              \norm{q}_{L^2(\Gamma)} &\le T_3 \norm{q}_{H^1(\Omega_p)} \quad \forall q\in X_p, \\
              \norm{\boldsymbol{\phi}}_{L^2(\Gamma)} &\le T_5 \norm{\boldsymbol{\phi}}_{H^1(\Omega_p)} \quad \forall \boldsymbol{\phi}\in\mathbf{H}_p.
          \end{align}
\end{itemize}

For the sake of simplicity, we define the following functions of time
\begin{equation}
\begin{split}
&\mathcal{C}^2_1(t):=\frac{1}{2}\norm{\mathbf{f}}_{L^2} ^2 + \frac{C_{P,p}^2}{2K_{\min}}\norm{g}_{L^2} ^2 + \frac{C_{P,p}^2}{2{a_p}_{min}}\norm{\phi}_{L^2} ^2,\\
&\mathcal{C}^2_2(t) := \tfrac12\|\mathbf{f}_t\|^2_{L^2}  + \frac{C_{P,p}^2}{2K_{\min}}\|g_t\|^2_{L^2}
+ \frac{C_{P,p}^2}{2{a_p}_{\min}}\|\phi_t\|^2_{L^2} ,\\
&\mathcal{C}_3^2:=\frac12\|\mathbf{f}(0)\|^2_{L^2}  + \frac{1}{2c_0}\|g(0)\|^2_{L^2} +\frac{1}{2c_0}\|\phi(0)\|^2_{L^2},
\end{split}
\label{eq:a30}
\end{equation}
small data condition holds:
\begin{equation}\label{cond:small}
\begin{aligned}
&
\|\mathcal{C}_1\|_{L^\infty(0,T)}^2+(1 + e^{T}+Te^{T})\|\mathcal{C}_1\|_{L^2(0,T)}^2  +  (1+Te^{T})  \|\mathcal{C}_2\|_{L^2(0,T)}^2 +  (1+Te^{T})\mathcal{C}_3^2<\mu_f\delta_1,\\
&
\|\mathcal{C}_1\|_{L^\infty(0,T)}^2+(1 + e^{T}+Te^{T})\|\mathcal{C}_1\|_{L^2(0,T)}^2  +  (1+Te^{T})  \|\mathcal{C}_2\|_{L^2(0,T)}^2 +  (1+Te^{T})\mathcal{C}_3^2<\frac{a_f}{2}\delta_2,
\end{aligned}
\end{equation}

\begin{align}
\delta_1 = \min\left\{
\frac{\mu_f}{4 C_{S,f}^2 C_K^3 },\;
\frac{a_f}{8C_{S,f}^2 C_K C_{P,f}},\;
\frac{a_f}{2C_{S,f}^2 C_K C_{P,f}} - \frac{2C_{11}^2 + 2C_{12}^2}{3C_{S,f}^2 C_K C_{P,f}\mu_f}
\right\},\
\delta_2 = \frac{\sqrt{2a_f \mu_f}}{4C_{S,f}^2 C_K C_{P,f}},
\label{eq:a31}
\end{align}
\begin{align*}
C_{11} := \frac{Gr}{Re^2} C_{S,f}^2 C_K, \qquad C_{12} := C_{S,f}^2 C_K.
\end{align*}

With these preliminary estimates and \eqref{cond:small} at hand, we are now in a position to establish the uniform a priori bounds for the Galerkin approximations.
\begin{theorem}
Suppose that $\mathbf{f} \in H^1(0,T;\mathbf{L}^2(\Omega_f))$, $g \in H^1(0,T;L^2(\Omega_p))$ and $\phi \in H^1(0,T;L^2(\Omega_p))$.
In addition, assume that the small data condition \eqref{cond:small} holds. Then, problem \eqref{eq:28}--\eqref{eq:33} has a unique solution $(\mathbf{u}_{f,m},\theta_{f,m}, \boldsymbol{\eta}_{p,m}, p_{p,m} , \theta_{p,m})$ in the interval $[0,T]$.
Furthermore, it satisfies the following bounds:
\begin{equation}
\begin{split}
&\frac{1}{2}\norm{\mathbf{u}_{f,m}}_{L^2}^2 + \frac{1}{2}\norm{\theta_{f,m}}_{L^2}^2 + \frac{1}{2}\norm{\bm{\xi}_{p,m}}_{L^2}^2 + \frac{c_0}{2}\norm{p_{p,m}}_{L^2}^2 + \frac{c_0}{2}\norm{\theta_{p,m}}_{L^2}^2  + \mu_p\norm{\mathbb{D}(\bm{\eta}_{p,m})}_{L^2}^2 \\
&\quad+ \frac{\lambda_p}{2}\norm{\nabla\cdot\bm{\eta}_{p,m}}_{L^2}^2 + \mu_f \norm{\mathbb{D}(\mathbf{u}_{f,m})}_{\Omegaf}^2 + \frac{a_f}{2}\norm{\nabla\theta_{f,m}}_{L^2}^2 \\
&\quad + \frac{\mathbf{K}}{2}\norm{\nabla p_{p,m}}_{L^2}^2 + \frac{\bm{a}_p}{2}\norm{\nabla\theta_{p,m}}_{L^2}^2 \leq ( 1+e^{T})\|\mathcal{C}_1\|_{L^2(0,T)}^2 ,
\end{split}
\label{eq:42}
\end{equation}

\begin{align}
\|\mathbb{D}(\mathbf{u}_{f,m})(t)\|_{L^2} < \delta_1,\qquad
\|\nabla\theta_{f,m}(t)\|_{L^2} <\delta_2 \label{eq:44} ,
\end{align}
and
\begin{equation}
\begin{split}
&\frac12\|\partial_t\mathbf{u}_{f,m}\|^2_{L^2}
+ \frac12\|\partial_t\theta_{f,m}\|^2_{L^2}
+ \frac12\|\partial_t\boldsymbol{\xi}_{p,m}\|^2_{L^2}
+ \frac{c_0}{2}\|\partial_t p_{p,m}\|^2_{L^2}
+ \frac{c_0}{2}\|\partial_t \theta_{p,m}\|^2_{L^2} + \mu_p\|\mathbb{D}(\boldsymbol{\xi}_{p,m})\|^2_{L^2}
\\
&\quad+ \frac{\lambda_p}{2}\|\nabla\!\cdot\!\boldsymbol{\xi}_{p,m}\|^2_{L^2}
+ \mu\|\mathbb{D}(\partial_t\mathbf{u}_{f,m})\|^2_{L^2(0,T;L^2(\Omf))} + \frac{a_f}{2}\|\nabla\partial_t\theta_{f,m}\|^2_{L^2(0,T;L^2(\Omf))}
\\
&\quad+ \frac{\mathbf{K}}{2}\|\nabla\partial_t p_{p,m}\|^2_{L^2(0,T;L^2(\Omp))}
+ \frac{\bm{a}_p}{2}\|\nabla\partial_t \theta_{p,m}\|^2_{L^2(0,T;L^2(\Omp))} \\
&\le (1+Te^{T}) \|\mathcal{C}_2\|_{L^2(0,T)}^2 + (1+Te^{T})\mathcal{C}_3^2. \label{eq:43}
\end{split}
\end{equation}
\end{theorem}
\begin{proof}\
\noindent \textbf{Proof of \eqref{eq:42}}\\
Take $\mathbf{v}_f = \mathbf{u}_{f,m}$, $w_f = \theta_{f,m}$, $\bm{\varphi}_p = \bm{\xi}_{p,m}$, $q_p = p_{p,m}$, and $w_p = \theta_{p,m}$ in \eqref{eq:28}--\eqref{eq:33}. After summation
we obtain:
\begin{align}
\frac{d}{dt}E_m(t) + D_m(t) + S_m(t) = R_m(t) - N_m(t),\label{eq:45}
\end{align}
where the $E_m(t)$, $D_m(t)$, $S_m(t)$, $ R_m(t)$, $N_m(t)$ are defined by
\begin{align*}
E_m(t)&=\bigl(\partial_t \mathbf{u}_{f,m},\mathbf{u}_{f,m}\bigr)_{\Omega_f}
+\bigl(\partial_t \theta_{f,m},\theta_{f,m}\bigr)_{\Omega_f}
+\bigl(\partial_t \boldsymbol{\xi}_{p,m},\boldsymbol{\xi}_{p,m}\bigr)_{\Omega_p}  +S_0\bigl(\partial_t p_{p,m},p_{p,m}\bigr)_{\Omega_p}
\notag\\
&\quad+a_0\bigl(\partial_t \theta_{p,m},\theta_{p,m}\bigr)_{\Omega_p}
-b_0\bigl(\partial_t \theta_{p,m},p_{p,m}\bigr)_{\Omega_p}
-b_0\bigl(\partial_t p_{p,m},\theta_{p,m}\bigr)_{\Omega_p}
\notag\\
&\quad
+2\mu_p \bigl(\mathbb{D}(\boldsymbol{\eta}_{p,m}),\mathbb{D}(\boldsymbol{\xi}_{p,m})\bigr)_{\Omega_p}
+\lambda_p \bigl(\nabla\cdot\boldsymbol{\eta}_{p,m},\nabla\cdot\boldsymbol{\xi}_{p,m}\bigr)_{\Omega_p},\\[4pt]
D_m(t)&=2\mu_f \|\mathbb{D}(\mathbf{u}_{f,m})\|^2_{\Omega_f}
+a_f \|\nabla\theta_{f,m}\|^2_{\Omega_f}  +\mathbf{K} \|\nabla p_{p,m}\|^2_{\Omega_p}
+\bm{a}_p \|\nabla\theta_{p,m}\|^2_{\Omega_p},\\[4pt]
S_m(t)&:=\frac{\alpha\sqrt{d}}{\sqrt{\operatorname{tr}(\mathbf{K})}}\Bigl[
\langle \mathbf{P}_{\boldsymbol{\tau}}(\mathbf{u}_{f,m}-\bm{\xi}_{p,m}),\mathbf{P}_{\boldsymbol{\tau}}(\mathbf{u}_{f,m})\rangle_\Gamma
- \langle \mathbf{P}_{\boldsymbol{\tau}}(\mathbf{u}_{f,m}-\bm{\xi}_{p,m}),\mathbf{P}_{\boldsymbol{\tau}}(\bm{\xi}_{p,m})\rangle_\Gamma \Bigr] \notag\\
&=\frac{\alpha\sqrt{d}}{\sqrt{\operatorname{tr}(\mathbf{K})}}\|\mathbf{P}_{\boldsymbol{\tau}}(\mathbf{u}_{f,m}-\bm{\xi}_{p,m})\|^2_\Gamma\ge 0,\\[4pt]
N_m(t) &:=\big((\mathbf{u}_{f,m}\cdot\nabla)\mathbf{u}_{f,m},\mathbf{u}_{f,m}\big)_{\Omega_f}
          + \big((\mathbf{u}_{f,m}\cdot\nabla)\theta_{f,m},\theta_{f,m}\big)_{\Omega_f},\\[4pt]
R_m(t) &:=\frac{Gr}{Re^2} (\bm{\xi}\theta_{f,m},\mathbf{u}_{f,m})_{\Omega_f}
+ (\mathbf{u}_{f,m}\cdot\bm{\xi},\theta_{f,m})_{\Omega_f}
+ (\mathbf{f},\bm{\xi}_{p,m})_{\Omega_p} + (g,p_{p,m})_{\Omega_p} + (\phi,\theta_{p,m})_{\Omega_p}.
\end{align*}
Next, we estimate $N_m(t)$ and $R_m(t)$ separately.\\
\noindent \textbf{Estimate of $N_m(t)$}\\
For the first term, using Sobolev inequality and Korn's inequality, we obtain
\begin{align*}
\big|((\mathbf{u}_{f,m}\cdot\nabla)\mathbf{u}_{f,m},\mathbf{u}_{f,m})\big|
\le \|\mathbf{u}_{f,m}\|_{L^4}\|\nabla\mathbf{u}_{f,m}\|_{L^2}\|\mathbf{u}_{f,m}\|_{L^4}
\le C_{S,f}^2 C_K^3 \|\D(\mathbf{u}_{f,m})\|_{L^2}^3.
\end{align*}

For the second term,
\begin{align*}
\big|((\mathbf{u}_{f,m}\cdot\nabla)\theta_{f,m},\theta_{f,m})\big|
\le \|\mathbf{u}_{f,m}\|_{L^4}\|\nabla\theta_{f,m}\|_{L^2}\|\theta_{f,m}\|_{L^4}
\le C_{S,f}^2 \|\mathbf{u}_{f,m}\|_{H^1} \|\nabla\theta_{f,m}\|_{L^2} \|\theta_{f,m}\|_{H^1},
\end{align*}
where
\begin{align*}
\|\theta_{f,m}\|_{H^1}^2
= \|\theta_{f,m}\|_{L^2}^2 + \|\nabla\theta_{f,m}\|_{L^2}^2 &\le C_P^2 \|\nabla\theta_{f,m}\|_{L^2}^2 + \|\nabla\theta_{f,m}\|_{L^2}^2 \le C_{P,f}^2 \|\nabla\theta_{f,m}\|_{L^2}^2 .
\end{align*}
Hence
\begin{align*}
\big|((\mathbf{u}_{f,m}\cdot\nabla)\theta_{f,m},\theta_{f,m})\big|
\le C_{S,f}^2 C_K C_{P,f} \|\D(\mathbf{u}_{f,m})\|_{L^2} \|\nabla\theta_{f,m}\|^2_{L^2}.
\end{align*}

Combining these estimates yields
\begin{align}
|N_m(t)| \le  C_{S,f}^2 C_K^3 \|\D(\mathbf{u}_{f,m})\|_{L^2}^3 + C_{S,f}^2 C_K C_{P,f} \|\D(\mathbf{u}_{f,m})\|_{L^2} \|\nabla\theta_{f,m}\|_{L^2}^2.
\end{align}
\noindent \textbf{Estimate of $\frac{Gr}{Re^2} (\bm{\xi}\theta_{f,m},\mathbf{u}_{f,m})_{\Omega_f}$ and $(\mathbf{u}_{f,m}\cdot\bm{\xi},\theta_{f,m})_{\Omega_f}$}\\
Using H\"older's inequality (with exponents $(2,4,4)$), Sobolev embedding, and Korn's inequality, we obtain
\begin{align*}
\Bigl|\frac{Gr}{Re^2} (\bm{\xi}\theta_{f,m},\mathbf{u}_{f,m})_{\Omega_f}\Bigr| &\le C_{11} \norm{\nabla\theta_{f,m}}_{L^2} \norm{\D(\mathbf{u}_{f,m})}_{L^2}, \\
\bigl| (\mathbf{u}_{f,m}\cdot\bm{\xi},\theta_{f,m})_{\Omega_f}\bigr| &\le C_{12} \norm{\nabla\theta_{f,m}}_{L^2} \norm{\D(\mathbf{u}_{f,m})}_{L^2} ,
\end{align*}
where
\begin{align*}
C_{11} := \frac{Gr}{Re^2} C_{S,f}^2 C_K, \qquad C_{12} := C_{S,f}^2 C_K.
\end{align*}
Applying Young's inequality $ab \le \frac{\epsilon}{2}a^2 + \frac{1}{2\epsilon}b^2$ and choosing
$\epsilon_1 = \frac{3\mu_f}{4C_{11}}$, $\epsilon_2 = \frac{3\mu_f}{4C_{12}}$ yields
\begin{align}
|\frac{Gr}{Re^2} (\bm{\xi}\theta_{f,m},\mathbf{u}_{f,m})_{\Omega_f}| + |(\mathbf{u}_{f,m}\cdot\bm{\xi},\theta_{f,m})_{\Omega_f}| \le \frac{3\mu_f}{4} \|\D(\mathbf{u}_{f,m})\|^2_{L^2} + \frac{2(C_{11}^2+C_{12}^2)}{3\mu_f} \|\nabla\theta_{f,m}\|^2_{L^2}.\label{eq:46}
\end{align}

\noindent \textbf{Estimate of $(\mathbf{f},\bm{\xi}_{p,m})_{\Omega_p}$}\\
\begin{align}
|(\mathbf{f},\bm{\xi}_{p,m})_{\Omega_p} | \le \norm{\mathbf{f}}_{L^2} \norm{\bm{\xi}_{p,m}}_{L^2}
      \le \frac12 \norm{\bm{\xi}_{p,m}}^2_{L^2} + \frac12 \norm{\mathbf{f}}^2_{L^2}. \label{eq:47}
\end{align}

\noindent \textbf{Estimate of $(g,p_{p,m})_{\Omega_p}$}\\
First use Poincaré inequality:
\begin{align*}
| (g,p_{p,m})_{\Omega_p} | \le \norm{g}_{L^2} \norm{p_{p,m}}_{L^2}
      \le C_{P,p} \norm{g}_{L^2} \norm{\nabla p_{p,m}}_{L^2}
      \le \frac{C_{P,p}}{\sqrt{K_{\min}}} \norm{g}_{L^2} \norm{\sqrt{\mathbf{K}}\nabla p_{p,m}}_{L^2}.
\end{align*}
Then Young's inequality yields
\begin{align}
|I_4| \le \frac{\mathbf{K}}{2} \norm{\nabla p_{p,m}}^2_{L^2} + \frac{C_{P,p}^2}{2K_{\min}} \norm{g}^2_{L^2}.\label{eq:48}
\end{align}

\noindent \textbf{Estimate of $(\phi,\theta_{p,m})_{\Omega_p}$}\\
Similarly,
\begin{align}
|(\phi,\theta_{p,m})_{\Omega_p}|
      \le \frac{\bm{a}_p}{2} \norm{\nabla\theta_{p,m}}^2_{L^2} + \frac{C_{P,p}^2}{2{a_p}_{min}} \norm{\phi}^2_{L^2} .\label{eq:49}\end{align}

Combining the estimates \eqref{eq:46}--\eqref{eq:49} in \eqref{eq:45}, we can get:
\begin{equation}
\begin{split}
&\frac{d}{dt}E_m(t) + D_m(t) + S_m(t)\\
&\le C_{S,f}^2 C_K^3 \|\D(\mathbf{u}_{f,m})\|_{L^2}^3 + C_{S,f}^2 C_K C_{P,f} \|\D(\mathbf{u}_{f,m})\|_{L^2} \|\nabla\theta_{f,m}\|_{L^2}^2+\frac12 \norm{\bm{\xi}_{p,m}}^2_{L^2} \\
&\quad+ \frac12 \norm{\mathbf{f}}^2_{L^2}
+\frac{\mathbf{K}}{2} \norm{\nabla p_{p,m}}^2_{L^2} + \frac{C_{P,p}^2}{2K_{\min}} \norm{g}^2_{L^2}+\frac{\bm{a}_p}{2} \norm{\nabla\theta_{p,m}}^2_{L^2} + \frac{C_{P,p}^2}{2{a_p}_{\min}} \norm{\phi}^2_{L^2}.
\end{split}
\label{eq:52}
\end{equation}
The only problematic terms on the right-hand side are
$\|\D(\mathbf{u}_{f,m})\|_{L^2}^3$, $\|\D(\mathbf{u}_{f,m})\|_{L^2}\|\nabla\theta_{f,m}\|_{L^2}^2$ and
$\|\bm{\xi}_{p,m}\|_{L^2}^2$. All other terms can be easily shifted to the left-hand side.
Observe that since $\mathbf{u}_{f,m}(0) = 0$ and $\mathbf{u}_{f,m}$ is continuous,
there exists a time $\bar{T}_m$, $0 < \bar{T}_m \le T_m$ such that
\begin{align}
\|\mathbb{D}(\mathbf{u}_{f,m})(t)\|_{L^2} <\delta_1
\quad \forall t\in[0,\bar{T}_m).\label{eq:5a4}
\end{align}
In fact, this condition holds true on $[0,{T}_m]$. For the sake of presentation, we postpone this proof to next. Using this condition, we obtain
\begin{equation}
\begin{aligned}
&\frac{d}{dt}E_m(t) + D_m(t) + S_m(t)\\
&\le \mu_f \|\D(\mathbf{u}_{f,m})\|_{L^2}^2 + \frac{a_f}{2}  \|\nabla\theta_{f,m}\|_{L^2}^2+\frac12 \norm{\bm{\xi}_{p,m}}^2_{L^2} + \frac12 \norm{\mathbf{f}}^2_{L^2}
+\frac{\mathbf{K}}{2} \norm{\nabla p_{p,m}}^2_{L^2}\\
&\quad + \frac{C_{P,p}^2}{2K_{\min}} \norm{g}^2_{L^2}+\frac{\bm{a}_p}{2} \norm{\nabla\theta_{p,m}}^2_{L^2} + \frac{C_{P,p}^2}{2{a_p}_{\min}} \norm{\phi}^2_{L^2}.
\end{aligned}
\label{eq:54}
\end{equation}

We note that:
\begin{gather}
(\partial_t \mathbf{u}_{f,m}, \mathbf{u}_{f,m}) = \frac12 \frac{d}{dt}\|\mathbf{u}_{f,m}\|^2_{L^2}, \
(\partial_t \theta_{f,m}, \theta_{f,m}) = \frac12 \frac{d}{dt}\|\theta_{f,m}\|^2_{L^2}, \
(\partial_t \bm{\xi}_{p,m}, \bm{\xi}_{p,m}) = \frac12 \frac{d}{dt}\|\bm{\xi}_{p,m}\|^2_{L^2}, \notag\\
S_0(\partial_t p_{p,m}, p_{p,m}) = \frac{S_0}{2}\frac{d}{dt}\|p_{p,m}\|^2_{L^2}, \
a_0(\partial_t \theta_{p,m}, \theta_{p,m}) = \frac{a_0}{2}\frac{d}{dt}\|\theta_{p,m}\|^2_{L^2}, \label{eq:060}
\end{gather}
\begin{align}
2\mu_p (\D(\bm{\eta}_{p,m}),\D(\bm{\xi}_{p,m})) + \lambda_p (\nabla\cdot\bm{\eta}_{p,m},\nabla\cdot\bm{\xi}_{p,m})
= \frac{d}{dt}\Bigl( \mu_p \|\D(\bm{\eta}_{p,m})\|^2_{L^2} + \frac{\lambda_p}{2} \|\nabla\cdot\bm{\eta}_{p,m}\|^2_{L^2} \Bigr).\label{eq:061}
\end{align}

By the thermodynamic stability condition $S_0 a_0 > b_0^2$,
the quadratic form $\frac{S_0}{2}\|p_{p,m}\|^2 + \frac{a_0}{2}\|\theta_{p,m}\|^2 - b_0(p_{p,m},\theta_{p,m})$
is positive definite, i.e., there exists a constant $c_0 > 0$ such that
\begin{align}
\frac{S_0}{2}\|p_{p,m}\|^2 + \frac{a_0}{2}\|\theta_{p,m}\|^2 - b_0(p_{p,m},\theta_{p,m}) \ge \frac{c_0}{2}\|p_{p,m}\|^2+\frac{c_0}{2}\|\theta_{p,m}\|^2.\label{eq:062}
\end{align}
For all $t \in [0, T_m]$, substituting \eqref{eq:060}--\eqref{eq:062} into \eqref{eq:54} yields
\begin{align}
&\frac{1}{2}\frac{d}{dt}\norm{\mathbf{u}_{f,m}}_{L^2}^2 + \frac{1}{2}\frac{d}{dt}\norm{\theta_{f,m}}_{L^2}^2 + \frac{1}{2}\frac{d}{dt}\norm{\bm{\xi}_{p,m}}_{L^2}^2 + \frac{c_0}{2}\frac{d}{dt}\norm{p_{p,m}}_{L^2}^2 + \frac{c_0}{2}\frac{d}{dt}\norm{\theta_{p,m}}_{L^2}^2 \notag\\
&\quad + \mu_p\frac{d}{dt}\norm{\mathbb{D}(\bm{\eta}_{p,m})}_{L^2}^2 + \frac{\lambda_p}{2}\frac{d}{dt}\norm{\nabla\cdot\bm{\eta}_{p,m}}_{L^2}^2 + \mu_f \norm{\mathbb{D}(\mathbf{u}_{f,m})}_{L^2}^2 + \frac{a_f}{2}\norm{\nabla\theta_{f,m}}_{L^2} \notag\\
&\quad + \frac{\mathbf{K}}{2}\norm{\nabla p_{p,m}}_{L^2}^2 + \frac{\bm{a}_p}{2}\norm{\nabla\theta_{p,m}}_{L^2}^2 \notag\\
&\leq \frac{1}{2}\norm{\bm{\xi}_{p,m}}_{L^2}^2 + \frac{1}{2}\norm{\mathbf{f}}_{L^2}^2 + \frac{(C_{P,p})^2}{2K_{\min}}\norm{g}_{L^2}^2 + \frac{(C_{P,p})^2}{2{a_p}_{min}}\norm{\phi}_{L^2}^2 \notag\\
&= \mathcal{C}_1^2(t) + \frac{1}{2}\norm{\bm{\xi}_{p,m}}_{L^2}^2.\label{eq:541}
\end{align}
Integrating \eqref{eq:541} over \((0,t)\) and using \eqref{eq:33}, then
retaining only the term involving \(\|\bm{\xi}_{p,m}\|_{L^2}^2\), we obtain
\begin{align*}
\norm{\bm{\xi}_{p,m}(t)}_{L^2}^2 \leq 2\|\mathcal{C}_1\|_{L^2(0,T)}^2 + \int_0^t \norm{\bm{\xi}_{p,m}(s)}_{L^2}^2 \,ds.
\end{align*}
Therefore, since $\bm{\eta}_{p,m} \in C^2(0,T; \bm{X}_p^m)$, $\norm{\bm{\xi}_{p,m}(t)}_{L^2}^2 \in C^1(0,T)$ and applying Lemma~\ref{thm:gronwall} with $\gamma(t) = \norm{\bm{\xi}_{p,m}(t)}_{L^2}^2$, $C = 1$ and $B = 2\|\mathcal{C}_1\|_{L^2(0,T)}^2$ yields
\begin{align*}
\norm{\bm{\xi}_{p,m}(t)}_{L^2}^2 \leq 2e^{T} \|\mathcal{C}_1\|_{L^2(0,T)}^2.
\end{align*}
Plugging this in \eqref{eq:541}, we have
\begin{align}
&\frac{1}{2}\frac{d}{dt}\norm{\mathbf{u}_{f,m}}_{L^2}^2 + \frac{1}{2}\frac{d}{dt}\norm{\theta_{f,m}}_{L^2}^2 + \frac{1}{2}\frac{d}{dt}\norm{\bm{\xi}_{p,m}}_{L^2}^2 + \frac{c_0}{2}\frac{d}{dt}\norm{p_{p,m}}_{L^2}^2 + \frac{c_0}{2}\frac{d}{dt}\norm{\theta_{p,m}}_{L^2}^2 \notag\\
&\quad + \mu_p\frac{d}{dt}\norm{\mathbb{D}(\bm{\eta}_{p,m})}_{L^2}^2 + \frac{\lambda_p}{2}\frac{d}{dt}\norm{\nabla\cdot\bm{\eta}_{p,m}}_{L^2}^2 + \mu_f \norm{\mathbb{D}(\mathbf{u}_{f,m})}_{L^2}^2 + \frac{a_f}{2}\norm{\nabla\theta_{f,m}}_{L^2}^2 \notag\\
&\quad + \frac{\mathbf{K}}{2}\norm{\nabla p_{p,m}}_{L^2}^2 + \frac{\bm{a}_p}{2}\norm{\nabla\theta_{p,m}}_{L^2}^2 \leq \mathcal{C}_1^2(t) + e^{T}\|\mathcal{C}_1\|_{L^2(0,T)}^2,\label{eq:56q}
\end{align}
for all $t \in [0, T_m]$.

After integrating with respect to $t$ and using \eqref{eq:33} implies the bound \eqref{eq:42}.

\noindent \textbf{Proof of \eqref{eq:43} and \eqref{eq:44}}\\
Recalling \eqref{eq:5a4}, assume for a contradiction that there exists
$T^* \in (0, T_m]$ such that for all $t \in [0, T^*)$, such that
\begin{align}
\label{eq:57}
\|\mathbf{D}(\mathbf{u}_{f,m})(t)\|_{L^2} < \delta_1,
\|\mathbf{D}(\mathbf{u}_{f,m})(T^*)\|_{L^2} =\delta_1,\quad\|\nabla\theta_{f,m}(t)\|^2_{L^2}<\delta_2,\|\nabla\theta_{f,m}(T^*)\|^2_{L^2}=\delta_2.
\end{align}

By differentiating each of the Galerkin equations \eqref{eq:28}--\eqref{eq:33}, we obtain:
\begin{align}
&(\partial_t\partial_t\mathbf{u}_{f,m},\mathbf{v}_f) + 2\mu_f(\mathbb{D}(\partial_t\mathbf{u}_{f,m}),\mathbb{D}(\mathbf{v}_f))
+ \langle\partial_t p_{p,m},\mathbf{v}_f\!\cdot\!\nff\rangle_\Gamma +((\partial_t\mathbf{u}_{f,m}\!\cdot\!\nabla)\mathbf{u}_{f,m},\mathbf{v}_f)\nonumber\\
&\quad +\frac{\alpha\sqrt{d}}{\sqrt{\tr(\mathbf{K})}}\langle\Ptau(\partial_t\mathbf{u}_{f,m}-\partial_t\boldsymbol{\xi}_{p,m}),
\Ptau(\mathbf{v}_f)\rangle_\Gamma
+((\mathbf{u}_{f,m}\!\cdot\!\nabla)\partial_t\mathbf{u}_{f,m},\mathbf{v}_f)  = \frac{Gr}{Re^2}(\boldsymbol{\xi}\partial_t\theta_{f,m},\mathbf{v}_f), \label{eq:D1}\\[4pt]
&(\partial_t\partial_t\theta_{f,m},w_f) + a_f(\nabla\partial_t\theta_{f,m},\nabla w_f)
+((\partial_t\mathbf{u}_{f,m}\!\cdot\!\nabla)\theta_{f,m},w_f)
+((\mathbf{u}_{f,m}\!\cdot\!\nabla)\partial_t\theta_{f,m},w_f) \nonumber\\
&\quad - a_f\langle\nabla\partial_t\theta_{f,m}\!\cdot\!\nff,w_f\rangle_\Gamma
= (\partial_t\mathbf{u}_{f,m}\!\cdot\!\boldsymbol{\xi},w_f), \label{eq:D2}\\[4pt]
&(\partial_t\partial_t\boldsymbol{\xi}_{p,m},\boldsymbol{\varphi}_p) + 2\mu_p(\mathbb{D}(\boldsymbol{\xi}_{p,m}),\mathbb{D}(\boldsymbol{\varphi}_p))
+ \lambda_p(\nabla\!\cdot\!\boldsymbol{\xi}_{p,m},\nabla\!\cdot\!\boldsymbol{\varphi}_p)  -\alpha(\nabla\!\cdot\!\boldsymbol{\varphi}_p,\partial_t p_{p,m}) -\beta(\nabla\!\cdot\!\boldsymbol{\varphi}_p,\partial_t \theta_{p,m})
\nonumber\\
&\quad- \langle\partial_t p_{p,m},\boldsymbol{\varphi}_p\!\cdot\!\nff\rangle_\Gamma -\frac{\alpha\sqrt{d}}{\sqrt{\tr(\mathbf{K})}}\langle\Ptau(\partial_t\mathbf{u}_{f,m}-\partial_t\boldsymbol{\xi}_{p,m}),
\Ptau(\boldsymbol{\varphi}_p)\rangle_\Gamma
= (\mathbf{f}_t,\boldsymbol{\varphi}_p), \label{eq:D3}\\[4pt]
&S_0(\partial_t\partial_t p_{p,m},q_p) - b_0(\partial_t\partial_t \theta_{p,m},q_p)
+ \alpha(\nabla\!\cdot\!\partial_t\boldsymbol{\xi}_{p,m},q_p) + \mathbf{K}(\nabla\partial_t p_{p,m},\nabla q_p) \nonumber\\
&\quad - \langle(\partial_t\mathbf{u}_{f,m}-\partial_t\boldsymbol{\xi}_{p,m})\!\cdot\!\nff,q_p\rangle_\Gamma
= (g_t,q_p), \label{eq:D4}\\[4pt]
&a_0(\partial_t\partial_t \theta_{p,m},w_p) - b_0(\partial_t\partial_t p_{p,m},w_p)
+ \beta(\nabla\!\cdot\!\partial_t\boldsymbol{\xi}_{p,m},w_p) + \bm{a}_p(\nabla\partial_t \theta_{p,m},\nabla w_p) \nonumber\\
&\quad - \bm{a}_p\langle\nabla\partial_t \theta_{p,m}\!\cdot\!\npp,w_p\rangle_\Gamma
= (\phi_t,w_p), \label{eq:D5}
\end{align}
where $\mathbf{f}_t=\partial_t\mathbf{f}$, $g_t=\partial_t g$, $\phi_t=\partial_t\phi$.

We select the test functions $\mathbf{v}_f = \partial_t\mathbf{u}_{f,m}$, $w_f = \partial_t\theta_{f,m}$,
$\boldsymbol{\varphi}_p = \partial_t\boldsymbol{\xi}_{p,m}$, $q_p = \partial_t p_{p,m}$, and $w_p = \partial_t \theta_{p,m}$;
substituting these into \eqref{eq:D1}--\eqref{eq:D5} and summing the five equations, we obtain:
\begin{equation}
\ddt\tilde{E}_m(t) + \tilde{D}_m(t) + \tilde{S}_m(t) = \tilde{R}_m(t) - \tilde{N}_m(t), \label{eq:61}
\end{equation}
where the $\tilde{E}_m(t)$ functional for the time derivatives is defined by
\begin{align*}
\tilde{E}_m(t) &:= \frac12\|\partial_t\mathbf{u}_{f,m}\|^2_{L^2}
+ \frac12\|\partial_t\theta_{f,m}\|^2_{L^2}
+ \frac12\|\partial_t\boldsymbol{\xi}_{p,m}\|^2_{L^2}
+ \frac{S_0}{2}\|\partial_t p_{p,m}\|^2_{L^2}
+ \frac{a_0}{2}\|\partial_t \theta_{p,m}\|^2_{L^2} \\
&\quad - b_0(\partial_t p_{p,m},\partial_t \theta_{p,m})_{\Omega_p}
+ \mu_p\|\mathbb{D}(\boldsymbol{\xi}_{p,m})\|^2_{L^2}
+ \frac{\lambda_p}{2}\|\nabla\!\cdot\!\boldsymbol{\xi}_{p,m}\|^2_{L^2}, \\
\tilde{D}_m(t) &:= 2\mu_f \|\mathbb{D}(\partial_t\mathbf{u}_{f,m})\|^2_{L^2} + a_f \|\nabla\partial_t\theta_{f,m}\|^2_{L^2}
+ \mathbf{K} \|\nabla \partial_t p_{p,m}\|^2_{L^2} + \bm{a}_p \|\nabla\partial_t \theta_{p,m}\|^2_{L^2}, \\
\tilde{S}_m(t) &:= \frac{\alpha\sqrt{d}}{\sqrt{\operatorname{tr}(\mathbf{K})}} \|\mathbf{P}_{\boldsymbol{\tau}}(\partial_t\mathbf{u}_{f,m}-\partial_t \bm{\xi}_{p,m})\|^2_\Gamma \ge 0, \\
\widetilde{R}_m(t) &= \frac{Gr}{Re^2} (\boldsymbol{\xi} \partial_t\theta_{f,m}, \partial_t\mathbf{u}_{f,m})
      + (\partial_t\mathbf{u}_{f,m}\!\cdot\!\boldsymbol{\xi}, \partial_t\theta_{f,m}) + (\mathbf{f}_t, \partial_t\boldsymbol{\xi}_{p,m}) + (g_t, \partial_t p_{p,m}) + (\phi_t, \partial_t \theta_{p,m}), \\
\widetilde{N}_m(t) &= (\partial_t\mathbf{u}_{f,m}\!\cdot\!\nabla \mathbf{u}_{f,m}, \partial_t\mathbf{u}_{f,m}) +(\mathbf{u}_{f,m}\!\cdot\!\nabla \partial_t\mathbf{u}_{f,m}, \partial_t\mathbf{u}_{f,m})\\&+ (\partial_t\mathbf{u}_{f,m}\!\cdot\!\nabla \theta_{f,m}, \partial_t\theta_{f,m}) + (\mathbf{u}_{f,m}\!\cdot\!\nabla \partial_t\theta_{f,m}, \partial_t\theta_{f,m}).
\end{align*}
Due to the thermo‑poroelastic stability condition $S_0a_0>b_0^2$, there exists a constant $c_0>0$
such that
\begin{equation}
\begin{aligned}
\label{eq:62ab}
\tilde{E}_m(t) &\ge \frac12\|\partial_t\mathbf{u}_{f,m}\|^2_{L^2} + \frac12\|\partial_t\theta_{f,m}\|^2_{L^2}
+ \frac12\|\partial_t\boldsymbol{\xi}_{p,m}\|^2_{L^2} + \frac{c_0}{2}\|\partial_t p_{p,m}\|^2_{L^2} \\
&\quad + \frac{c_0}{2}\|\partial_t \theta_{p,m}\|^2_{L^2}
+ \mu_p\|\mathbb{D}(\boldsymbol{\xi}_{p,m})\|^2_{L^2} + \frac{\lambda_p}{2}\|\nabla\!\cdot\!\boldsymbol{\xi}_{p,m}\|^2_{L^2}.
\end{aligned}
\end{equation}
\noindent \textbf{Estimate of \(\tilde{N}_m(t)\)}
Using H\"older's inequality, Sobolev embedding, Korn's inequality, and the Poincar\'e inequality, we obtain
\begin{equation}
\begin{aligned}
\bigl|((\partial_t\mathbf{u}_{f,m}\!\cdot\!\nabla)\mathbf{u}_{f,m},\partial_t\mathbf{u}_{f,m})\bigr|
&\le C_{S,f}^2 C_K^3\,\delta_1\,\|\mathbb{D}(\partial_t\mathbf{u}_{f,m})\|^2_{L^2}, \\
\bigl|((\mathbf{u}_{f,m}\!\cdot\!\nabla)\partial_t\mathbf{u}_{f,m},\partial_t\mathbf{u}_{f,m})\bigr|
&\le C_{S,f}^2 C_K^3\,\delta_1\,\|\mathbb{D}(\partial_t\mathbf{u}_{f,m})\|^2_{L^2}, \\
\bigl|((\partial_t\mathbf{u}_{f,m}\!\cdot\!\nabla)\theta_{f,m},\partial_t\theta_{f,m})\bigr|
&\le C_{S,f}^2 C_K C_{P,f}\,\delta_2\,\|\mathbb{D}(\partial_t\mathbf{u}_{f,m})\|_{L^2}\,\|\nabla\partial_t\theta_{f,m}\|_{L^2}, \\
\bigl|((\mathbf{u}_{f,m}\!\cdot\!\nabla)\partial_t\theta_{f,m},\partial_t\theta_{f,m})\bigr|
&\le C_{S,f}^2 C_K C_{P,f}\,\delta_1\,\|\nabla\partial_t\theta_{f,m}\|^2_{L^2},
\end{aligned}
\label{eq:N_estimates}
\end{equation}
summing the four estimates yields
\begin{align*}
|\tilde{N}_m(t)| &\le 2C_{S,f}^2 C_K^3\delta_1 \|\mathbb{D}(\partial_t\mathbf{u}_{f,m})\|^2_{L^2}
+ C_{S,f}^2 C_K C_{P,f}\delta_2 \|\mathbb{D}(\partial_t\mathbf{u}_{f,m})\|_{L^2}\,\|\nabla\partial_t\theta_{f,m}\|_{L^2}
\\&+ C_{S,f}^2 C_K C_{P,f}\delta_1 \|\nabla\partial_t\theta_{f,m}\|^2_{L^2}.
\end{align*}
Applying Young's inequality to the mixed term, with a weight chosen to facilitate later absorption into the left‑hand side, we obtain
\begin{align}
|\tilde{N}_m(t)| &\le \bigl(2C_{S,f}^2 C_K^3\delta_1 + \tfrac{\mu_f}{4}\bigr)\|\mathbb{D}(\partial_t\mathbf{u}_{f,m})\|^2_{L^2}
+ \bigl(C_{S,f}^2 C_K C_{P,f}\delta_1 + \frac{(C_{S,f}^2 C_K C_{P,f}\delta_2)^2}{\mu_f}\bigr)\|\nabla\partial_t\theta_{f,m}\|^2_{L^2} \notag\\
&\le \frac{3\mu_f}{4}\|\mathbb{D}(\partial_t\mathbf{u}_{f,m})\|^2_{L^2}
+ \frac{a_f}{4} \|\nabla\partial_t\theta_{f,m}\|^2_{L^2}. \label{eq:Nest}
\end{align}
\noindent \textbf{Estimate of $\tilde{R}_m(t)$}
Applying H\"older's inequality, Korn's inequality, Poincar\'e's inequality, and Young's inequality to each term, we obtain the following estimate:
\begin{equation}
\begin{split}
|\tilde{R}_m(t)| &\le \frac{\epsilon_1+\epsilon_2}{2}\|\mathbb{D}(\partial_t\mathbf{u}_{f,m})\|^2_{L^2}+\Bigl(\frac{C_{11}^2}{2\epsilon_1}+\frac{C_{12}^2}{2\epsilon_2}\Bigr)\|\nabla\partial_t\theta_{f,m}\|^2_{L^2} \\
&\quad + \frac12\|\partial_t\boldsymbol{\xi}_{p,m}\|^2_{L^2} + \frac{\mathbf{K}}{2}\|\nabla\partial_t p_{p,m}\|^2_{L^2}
+ \frac{\bm{a}_p}{2}\|\nabla\partial_t \theta_{p,m}\|^2_{L^2} + \mathcal{C}_2^2(t), \label{eq:Rmest}
\end{split}
\end{equation}
where
\begin{align*}
\mathcal{C}_2^2(t) := \frac12\|\mathbf{f}_t\|^2_{L^2} + \frac{C_{P,p}^2}{2K_{\min}}\|g_t\|^2_{L^2}
+ \frac{C_{P,p}^2}{{2a_p}_{min}}\|\phi_t\|^2_{L^2}.
\end{align*}

Insert the estimates \eqref{eq:62ab}, \eqref{eq:Nest} and \eqref{eq:Rmest} in \eqref{eq:61} and note that $\tilde{S}_m(t)\ge0$. We obtain
\begin{align*}
&\frac12\frac{d}{dt}\|\partial_t\mathbf{u}_{f,m}\|^2_{L^2} + \frac12\frac{d}{dt}\|\partial_t\theta_{f,m}\|^2_{L^2}
+ \frac12\frac{d}{dt}\|\partial_t\boldsymbol{\xi}_{p,m}\|^2_{L^2} + \frac{c_0}{2}\frac{d}{dt}\|\partial_t p_{p,m}\|^2_{L^2}
+ \frac{c_0}{2}\frac{d}{dt}\|\partial_t \theta_{p,m}\|^2_{L^2}
\\
&\quad+ \mu_p\frac{d}{dt}\|\mathbb{D}(\boldsymbol{\xi}_{p,m})\|^2_{L^2}  + \frac{\lambda_p}{2}\frac{d}{dt}\|\nabla\!\cdot\!\boldsymbol{\xi}_{p,m}\|^2_{L^2}+2\mu_f \|\mathbb{D}(\partial_t\mathbf{u}_{f,m})\|^2_{L^2} + a_f \|\nabla\partial_t\theta_{f,m}\|^2_{L^2}
\\
&\quad+ \mathbf{K} \|\nabla \partial_t p_{p,m}\|^2_{L^2} + \bm{a}_p \|\nabla\partial_t \theta_{p,m}\|^2_{L^2}\\
&\le \frac{\epsilon_1+\epsilon_2}{2}\|\mathbb{D}(\partial_t\mathbf{u}_{f,m})\|^2_{L^2}+\Bigl(\frac{C_{11}^2}{2\epsilon_1}+\frac{C_{12}^2}{2\epsilon_2}\Bigr)\|\nabla\partial_t\theta_{f,m}\|^2_{L^2}
+ \tfrac12\|\partial_t\boldsymbol{\xi}_{p,m}\|^2_{L^2} + \tfrac{\mathbf{K}}{2}\|\nabla\partial_t p_{p,m}\|^2_{L^2}
\\
&\quad+ \frac{\bm{a}_p}{2}\|\nabla\partial_t \theta_{p,m}\|^2_{L^2} + \mathcal{C}_2^2(t)  +  \frac{3\mu_f}{4}\|\mathbb{D}(\partial_t\mathbf{u}_{f,m})\|^2_{L^2}
+ \frac{a_f}{4}\|\nabla\partial_t\theta_{f,m}\|^2_{L^2}.
\end{align*}
We choose the free parameters
\begin{align*}
\epsilon_1 = \frac{2C_{11}^2}{a_f},\qquad \epsilon_2 = \frac{2C_{12}^2}{a_f},
\end{align*}
which makes \(\frac{C_{11}^2}{2\epsilon_1}+\frac{C_{12}^2}{2\epsilon_2} = \frac{a_f}{4}\).
we obtain the simplified inequality
\begin{equation}
\begin{split}
&\frac12\frac{d}{dt}\|\partial_t\mathbf{u}_{f,m}\|^2_{L^2} + \frac12\frac{d}{dt}\|\partial_t\theta_{f,m}\|^2_{L^2}
+ \frac12\frac{d}{dt}\|\partial_t\boldsymbol{\xi}_{p,m}\|^2_{L^2} + \frac{c_0}{2}\frac{d}{dt}\|\partial_t p_{p,m}\|^2_{L^2}
+ \frac{c_0}{2}\frac{d}{dt}\|\partial_t \theta_{p,m}\|^2_{L^2}
\\
&\quad+ \mu_p\frac{d}{dt}\|\mathbb{D}(\boldsymbol{\xi}_{p,m})\|^2_{L^2}  + \frac{\lambda_p}{2}\frac{d}{dt}\|\nabla\!\cdot\!\boldsymbol{\xi}_{p,m}\|^2 _{L^2} +\mu_f\|\mathbb{D}(\partial_t\mathbf{u}_{f,m})\|^2_{L^2} + \frac{a_f}{2}\|\nabla\partial_t\theta_{f,m}\|^2_{L^2} \\
&\quad+ \frac{\mathbf{K}}{2}\|\nabla\partial_t p_{p,m}\|^2_{L^2} + \frac{\bm{a}_p}{2}\|\nabla\partial_t \theta_{p,m}\|^2_{L^2}
\le \frac12\|\partial_t\boldsymbol{\xi}_{p,m}\|^2_{L^2} + \mathcal{C}_2^2(t),
\end{split}
\label{eq:66}
\end{equation}

Next we need a bound for $\partial_t\boldsymbol{\xi}_{p,m}$. Since by assumption $\partial_t\boldsymbol{\xi}_{p,m}$ is continuous,
we can use Gronwall's inequality again. Integrating from $0$ to $t$ for any $t \in [0, T^*]$, we obtain
\begin{equation}
\begin{split}
&\frac12\|\partial_t\mathbf{u}_{f,m}\|^2_{L^2} + \frac12\|\partial_t\theta_{f,m}\|^2_{L^2}
+ \frac12\|\partial_t\boldsymbol{\xi}_{p,m}\|^2_{L^2} + \frac{c_0}{2}\|\partial_t p_{p,m}\|^2_{L^2}
+ \frac{c_0}{2}\|\partial_t \theta_{p,m}\|^2_{L^2}
\\
&\quad+ \mu_p\|\mathbb{D}(\boldsymbol{\xi}_{p,m})\|^2_{L^2}  + \frac{\lambda_p}{2}\|\nabla\!\cdot\!\boldsymbol{\xi}_{p,m}\|^2_{L^2}
\\&\le\frac12\|\partial_t\mathbf{u}_{f,m}(0)\|^2_{L^2} + \frac12\|\partial_t\theta_{f,m}(0)\|^2_{L^2}
+ \frac12\|\partial_t\boldsymbol{\xi}_{p,m}(0)\|^2_{L^2} + \frac{c_0}{2}\|\partial_t p_{p,m}(0)\|^2_{L^2}
 \\
&\quad + \frac{c_0}{2}\|\partial_t \theta_{p,m}(0)\|^2_{L^2}
+\|\mathcal{C}_2\|^2_{L^2(0,T)}+ \frac{1}{2}\int_0^t \norm{\partial_t\boldsymbol{\xi}_{p,m}\|^2(s)}_{L^2(\Omp)}^2 \,ds,
\end{split}
\label{eq:67}
\end{equation}
We derive bounds for the time derivatives of the Galerkin approximations at $t=0$. These estimates are needed later to apply Gronwall's inequality to the time derivatives.\\
\noindent \textbf{Estimate for $\partial_t^2 {\boldsymbol{\eta}}_{p,m}(0)$}\\
Choosing $\bm{\varphi}_p = \partial_t\boldsymbol{\xi}_{p,m}(0)=\partial_t\boldsymbol{\xi}_{p,m}(0)$ in \eqref{eq:30} and using the fact that all terms without time derivatives vanish at $t=0$, we obtain
\[
\bigl(\partial_t\boldsymbol{\xi}_{p,m}(0),\partial_t\boldsymbol{\xi}_{p,m}(0)\bigr)_{\Omega_p}
= \bigl(\mathbf{f}(0),\partial_t\boldsymbol{\xi}_{p,m}(0)\bigr)_{\Omega_p}.
\]
Since $\partial_t\boldsymbol{\xi}_{p,m}(0)=\partial_t\boldsymbol{\xi}_{p,m}(0)$, the left-hand side equals $\|\partial_t\boldsymbol{\xi}_{p,m}(0)\|_{L^2}^2$. By the Cauchy–Schwarz inequality,
\begin{align*}
\|\partial_t\boldsymbol{\xi}_{p,m}(0)\|_{L^2}^2 \le \|\mathbf{f}(0)\|_{L^2}\,\|\partial_t\boldsymbol{\xi}_{p,m}(0)\|_{L^2},
\end{align*}
hence
\begin{align}
\|\partial_t\boldsymbol{\xi}_{p,m}(0)\|_{L^2} \le \|\mathbf{f}(0)\|_{L^2}.\label{eq:a1}
\end{align}

\noindent \textbf{Estimates for $\partial_t p_{p,m}(0)$ and $\partial_t\theta_{p,m}(0)$}\\

We use the pressure equation \eqref{eq:31} and the solid temperature equation \eqref{eq:32}. At $t=0$ we get
\begin{align}
S_0\bigl(\partial_t p_{p,m}(0),q_p\bigr)_{\Omega_p} - b_0\bigl(\partial_t\theta_{p,m}(0),q_p\bigr)_{\Omega_p}
&= \bigl(g(0),q_p\bigr)_{\Omega_p}, \label{pp0} \\
a_0\bigl(\partial_t\theta_{p,m}(0),w_p\bigr)_{\Omega_p} - b_0\bigl(\partial_t p_{p,m}(0),w_p\bigr)_{\Omega_p}
&= \bigl(\phi(0),w_p\bigr)_{\Omega_p}. \label{tp0}
\end{align}
Take $q_p = \partial_t p_{p,m}(0)$ in \eqref{pp0} and $w_p = \partial_t\theta_{p,m}(0)$ in \eqref{tp0}. Adding the two equalities yields
\begin{align*}
&S_0\|\partial_t p_{p,m}(0)\|_{L^2}^2 + a_0\|\partial_t\theta_{p,m}(0)\|_{L^2}^2
-2b_0\bigl(\partial_t p_{p,m}(0),\partial_t\theta_{p,m}(0)\bigr)_{\Omega_p}\\&
= \bigl(g(0),\partial_t p_{p,m}(0)\bigr)_{\Omega_p}
+ \bigl(\phi(0),\partial_t\theta_{p,m}(0)\bigr)_{\Omega_p},
\end{align*}
and
\begin{align*}
&S_0\|\partial_t p_{p,m}(0)\|_{L^2}^2 + a_0\|\partial_t\theta_{p,m}(0)\|_{L^2}^2
-2b_0\bigl(\partial_t p_{p,m}(0),\partial_t\theta_{p,m}(0)\bigr)_{\Omega_p}\\& \ge c_0\bigl(\|\partial_t p_{p,m}(0)\|^2_{L^2}+\|\partial_t\theta_{p,m}(0)\|^2_{L^2}\bigr).
\end{align*}
Hence
\begin{align*}
c_0\bigl(\|\partial_t p_{p,m}(0)\|^2_{L^2}+\|\partial_t\theta_{p,m}(0)\|^2_{L^2}\bigr) \le \|g(0)\|_{L^2}\,\|\partial_t p_{p,m}(0)\|_{L^2} + \|\phi(0)\|_{L^2}\,\|\partial_t\theta_{p,m}(0)\|_{L^2}.
\end{align*}
Applying Young's inequality, $\|g(0)\|_{L^2}\|\partial_t p_{p,m}(0)\|_{L^2}\le \frac{c_0}{2}\|\partial_t p_{p,m}(0)\|^2_{L^2}+\frac{1}{2c_0}\|g(0)\|^2_{L^2}$, and similarly for the second term, we obtain
\begin{align}
\frac{c_0}{2}\bigl(\|\partial_t p_{p,m}(0)\|^2_{L^2}+\|\partial_t\theta_{p,m}(0)\|^2_{L^2}\bigr) \le \frac{1}{2c_0}\bigl(\|g(0)\|^2_{L^2}+\|\phi(0)\|^2_{L^2}\bigr).\label{eq:a2}
\end{align}
\noindent \textbf{Estimates for $\partial_t\mathbf{u}_{f,m}(0)$ and $\partial_t\theta_{f,m}(0)$}\\

The fluid momentum equation \eqref{eq:28} at $t=0$, with test function $\mathbf{v}_f=\partial_t\mathbf{u}_{f,m}(0)$, reduces to
\begin{align*}
\|\partial_t\mathbf{u}_{f,m}(0)\|_{L^2}^2 = \frac{G r}{Re^2}\bigl(\boldsymbol{\xi}\,\theta_{f,m}(0),\partial_t\mathbf{u}_{f,m}(0)\bigr)_{\Omega_f}.
\end{align*}
Since $\theta_{f,m}(0)=0$, the right‑hand side vanishes, giving
\begin{align}
\|\partial_t\mathbf{u}_{f,m}(0)\|_{L^2}=0. \label{eq:a3}
\end{align}
Similarly, the fluid temperature equation \eqref{eq:29} with $w_f=\partial_t\theta_{f,m}(0)$ and the initial condition $\theta_{f,m}(0)=0$, $\mathbf{u}_{f,m}(0)=0$ yields
\begin{align}
\|\partial_t\theta_{f,m}(0)\|_{L^2}^2 =0.\label{eq:a4}
\end{align}

Using \eqref{eq:a1}, \eqref{eq:a2}, \eqref{eq:a3} and \eqref{eq:a4} in\eqref{eq:67} we get
\begin{equation}
\begin{split}
&\frac12\|\partial_t\mathbf{u}_{f,m}\|^2_{L^2}+ \frac12\|\partial_t\theta_{f,m}\|^2_{L^2}
+ \frac12\|\partial_t\boldsymbol{\xi}_{p,m}\|^2_{L^2} + \frac{c_0}{2}\|\partial_t p_{p,m}\|^2_{L^2}
+ \frac{c_0}{2}\|\partial_t \theta_{p,m}\|^2_{L^2}
+ \mu_p\|\mathbb{D}(\boldsymbol{\xi}_{p,m})\|^2_{L^2} \\
&\quad + \frac{\lambda_p}{2}\|\nabla\!\cdot\!\boldsymbol{\xi}_{p,m}\|^2_{L^2}
\\&\le\frac12\|\mathbf{f}(0)\|^2_{L^2} + \frac{1}{2c_0}\|g(0)\|^2_{L^2}+\frac{1}{2c_0}\|\phi(0)\|^2_{L^2}
+\|\mathcal{C}_2\|^2_{L^2(0,T)}+ \frac{1}{2}\int_0^t \|\partial_t\boldsymbol{\xi}_{p,m}(s)\|_{L^2}^2 \,ds
\\&=\mathcal{C}_3^2
+\|\mathcal{C}_2\|^2_{L^2(0,T)}+ \frac{1}{2}\int_0^t \|\partial_t\boldsymbol{\xi}_{p,m}(s)\|_{L^2}^2 \,ds,
\end{split}
\end{equation}
where
\begin{align*}
\mathcal{C}_3^2=\frac12\|\mathbf{f}(0)\|^2_{L^2} + \frac{1}{2c_0}\|g(0)\|^2_{L^2}+\frac{1}{2c_0}\|\phi(0)\|^2_{L^2}.
\end{align*}
Neglecting all terms on the left-hand side of the above equation except the second one, we obtain
\begin{align*}
\|\partial_t\boldsymbol{\xi}_{p,m}\|^2_{L^2}
\leq2\mathcal{C}_3^2
+2\|\mathcal{C}_2\|^2_{L^2(0,T)}+ \int_0^t \|\partial_t\boldsymbol{\xi}_{p,m}\|_{L^2}^2(s)\, ds,
\end{align*}
Using Gronwall's inequality \cref{thm:gronwall} with
$B = 2 \mathcal{C}_3^2 + 2 \|\mathcal{C}_2\|_{L^2(0,T)}^2$, $C = 1$ and $\gamma(t) = \|\partial_t\boldsymbol{\xi}_{p,m}(t)\|_{L^2}^2$ we obtain
\begin{align}
\|\partial_t\boldsymbol{\xi}_{p,m}(t)\|_{L^2}^2 \le e^{T} \left( 2\|\mathcal{C}_2\|_{L^2(0,T)}^2 + 2 \mathcal{C}_2^3 \right).
\end{align}
Therefore, putting this in \eqref{eq:66} we have
\begin{equation}
\begin{split}
&\frac12\frac{d}{dt}\|\partial_t\mathbf{u}_{f,m}\|^2_{L^2}
+ \frac12\frac{d}{dt}\|\partial_t\theta_{f,m}\|^2_{L^2}+ \frac12\frac{d}{dt}\|\partial_t\boldsymbol{\xi}_{p,m}\|^2_{L^2}
+ \frac{c_0}{2}\frac{d}{dt}\|\partial_t p_{p,m}\|^2_{L^2}
 \\
&\quad+ \frac{c_0}{2}\frac{d}{dt}\|\partial_t \theta_{p,m}\|^2_{L^2} + \mu_p\frac{d}{dt}\|\mathbb{D}(\boldsymbol{\xi}_{p,m})\|^2_{L^2}
+ \frac{\lambda_p}{2}\frac{d}{dt}\|\nabla\!\cdot\!\boldsymbol{\xi}_{p,m}\|^2_{L^2} \\
&\quad+ \mu\|\mathbb{D}(\partial_t\mathbf{u}_{f,m})\|^2_{L^2}
+ \frac{a_f}{2}\|\nabla\partial_t\theta_{f,m}\|^2_{L^2}+ \frac{\mathbf{K}}{2}\|\nabla\partial_t p_{p,m}\|^2_{L^2}
+ \frac{\bm{a}_p}{2}\|\nabla\partial_t \theta_{p,m}\|^2_{L^2} \\
&\le e^{T} \left( \|\mathcal{C}_2\|_{L^2(0,T)}^2 + \mathcal{C}_3^2 \right) + \mathcal{C}_2^2(t), \label{eq:77}
\end{split}
\end{equation}
Integrating \eqref{eq:77} from 0 to t in $[0, T^*]$ and recalling the bounds \eqref{eq:a1}, \eqref{eq:a2}, \eqref{eq:a3} and \eqref{eq:a4} we obtain
\begin{equation}
\begin{split}
&\frac12\|\partial_t\mathbf{u}_{f,m}\|^2_{L^2} + \frac12\|\partial_t\theta_{f,m}\|^2_{L^2}
+ \frac12\|\partial_t\boldsymbol{\xi}_{p,m}\|^2_{L^2} + \frac{c_0}{2}\|\partial_t p_{p,m}\|^2_{L^2}
 + \frac{c_0}{2}\|\partial_t \theta_{p,m}\|^2_{L^2}
\\
&\quad+ \mu_p\|\mathbb{D}(\boldsymbol{\xi}_{p,m})\|^2_{L^2} + \frac{\lambda_p}{2}\|\nabla\!\cdot\!\boldsymbol{\xi}_{p,m}\|^2_{L^2}+\int_0^t\mu\|\mathbb{D}(\partial_t\mathbf{u}_{f,m})\|^2_{L^2}dt \\
&\quad+ \int_0^t\frac{a_f}{2}\|\nabla\partial_t\theta_{f,m}\|^2_{L^2}dt + \int_0^t\frac{\mathbf{K}}{2}\|\nabla\partial_t p_{p,m}\|^2_{L^2}dt + \int_0^t\frac{\bm{a}_p}{2}\|\nabla\partial_t \theta_{p,m}\|^2_{L^2}dt
\\&\le
\frac12\|\mathbf{f}_0\|^2_{L^2}+ \frac{1}{2c_0}\|g(0)\|^2_{L^2}+\frac{1}{2c_0}\|\phi(0)\|^2_{L^2}
+Te^{T} \left( \|\mathcal{C}_2\|_{L^2(0,T)}^2 +  \mathcal{C}_3^2 \right) + \|\mathcal{C}_2\|_{L^2(0,T)}^2
\\&\le
(1+Te^{T})  \|\mathcal{C}_2\|_{L^2(0,T)}^2 +  (1+Te^{T})\mathcal{C}_3^2   .\label{eq:81}
\end{split}
\end{equation}

Using Young’s inequality, \eqref{eq:54} yields
\begin{equation}
\begin{aligned}
&\mu_f \|\mathbb{D}(\mathbf{u}_{f,m})\|^2_{L^2}
+\frac{a_f}{2}\|\nabla\theta_{f,m}\|^2_{L^2}  \\
&\le \frac12\|\mathbf{u}_{f,m}\|^2_{L^2} + \frac12\|\theta_{f,m}\|^2_{L^2}
+ \frac{1}{2}\|\bm{\xi}_{p,m}\|^2_{L^2} + \frac{c_0}{2}\|p_{p,m}\|^2_{L^2}
+ \frac{c_0}{2}\|\theta_{p,m}\|^2_{L^2}  + \mu_p \|\mathbb{D}(\bm{\eta}_{p,m})\|^2_{L^2} \\
&\quad + \frac{\lambda_p}{2}\|\nabla\cdot\bm{\eta}_{p,m}\|^2_{L^2} +\frac12\|\partial_t\mathbf{u}_{f,m}\|^2_{L^2} + \frac12\|\partial_t\theta_{f,m}\|^2_{L^2}
+ \frac{1}{2}\|\partial_t \bm{\xi}_{p,m}\|^2_{L^2} + \frac{c_0}{2}\|\partial_t p_{p,m}\|^2_{L^2}
+ \frac{c_0}{2}\|\partial_t \theta_{p,m}\|^2_{L^2}  \\
&\quad + \mu_p \|\mathbb{D}(\partial_t \bm{\eta}_{p,m})\|^2_{L^2} + \frac{\lambda_p}{2}\|\nabla\cdot\partial_t \bm{\eta}_{p,m}\|^2_{L^2} +\mathcal{C}_1^2(t) + e^{T}\|\mathcal{C}_1\|_{L^2(0,T)}^2 ,\label{eq:82}
\end{aligned}
\end{equation}
we recall \eqref{eq:52} and use assumption \eqref{eq:57} to get \eqref{eq:56q} for all $t \in [0, T^*]$. Neglecting the terms other than the ones with $\frac{d}{dt}$ on the left hand side, we integrate \eqref{eq:56q} from $0$ to $t$ where $t \in [0, T^*]$. Since the initial conditions are zero, we have
\begin{align}
&\frac{1}{2}\norm{\mathbf{u}_{f,m}}_{L^2}^2 + \frac{1}{2}\norm{\theta_{f,m}}_{L^2}^2 + \frac{1}{2}\norm{\bm{\xi}_{p,m}}_{L^2}^2 + \frac{c_0}{2}\norm{p_{p,m}}_{L^2}^2 + \frac{c_0}{2}\norm{\theta_{p,m}}_{L^2}^2 \notag\\
&\quad + \mu_p\norm{D(\bm{\eta}_{p,m})}_{L^2}^2 + \frac{\lambda_p}{2}\norm{\nabla\cdot\bm{\eta}_{p,m}}_{L^2}^2 \leq (1 + Te^{T})\|\mathcal{C}_1\|_{L^2(0,T)}^2 .\label{eq:83}
\end{align}
Therefore, using \eqref{eq:81} and \eqref{eq:83} in \eqref{eq:82}, we get
\begin{equation}
\begin{aligned}
\|\mathbb{D}(\mathbf{u}_{f,m})\|^2_{L^2}&\le
\frac{1}{\mu_f}\Bigl(\mathcal{C}_1^2(t)+(1 + e^{T}+Te^{T})\|\mathcal{C}_1\|_{L^2(0,T)}^2 \\
&\qquad + (1+Te^{T}) \|\mathcal{C}_2\|_{L^2(0,T)}^2 + (1+Te^{T})\mathcal{C}_3^2\Bigr), \\[2mm]
\|\nabla\theta_{f,m}\|^2_{L^2}&\le
\frac{2}{a_f}\Bigl(\mathcal{C}_1^2(t)+(1 + e^{T}+Te^{T})\|\mathcal{C}_1\|_{L^2(0,T)}^2 \\
&\qquad + (1+Te^{T}) \|\mathcal{C}_2\|_{L^2(0,T)}^2 + (1+Te^{T})\mathcal{C}_3^2\Bigr),
\end{aligned}
\end{equation}
for all $t \in [0, T^*]$. Therefore, if the data satisfy
\begin{align*}
&
\frac{1}{\mu_f}\left(\mathcal{C}_1^2(t)+(1 + e^{T}+Te^{T})\|\mathcal{C}_1\|_{L^2(0,T)}^2  +  (1+Te^{T})  \|\mathcal{C}_2\|_{L^2(0,T)}^2 +  (1+Te^{T})\mathcal{C}_3^2\right)<\delta_1^2,\notag\\
&
\frac{2}{a_f}\left(\mathcal{C}_1^2(t)+(1 + e^{T}+Te^{T})\|\mathcal{C}_1\|_{L^2(0,T)}^2  +  (1+Te^{T})  \|\mathcal{C}_2\|_{L^2(0,T)}^2 +  (1+Te^{T})\mathcal{C}_3^2\right)<\delta_2^2,
\end{align*}
then we have
\begin{align}
\|\mathbb{D}(\mathbf{u}_{f,m})\|_{L^2}<\delta_1,\quad
\|\nabla\theta_{f,m}\|_{L^2}<\delta_2,
\end{align}
for all $t \in [0, T^*]$ which contradicts with the assumption.
\end{proof}
\subsection{Passing to the limit}\label{subsec:limit}

The uniform a priori estimates established in the previous section allow us to pass to the limit
\(m\to\infty\) and obtain a solution to the original weak formulation. The key idea is to exploit the boundedness of the sequences in appropriate Sobolev spaces to extract weakly convergent subsequences. To handle the nonlinear convective terms, we require stronger convergence, which will be obtained via a compactness argument (Lemma \ref{lem:compactness}) that relies on the uniform bounds on the time derivatives.

Since we established the boundedness of the Galerkin solution $(\mathbf{u}_{f,m},\theta_{f,m},\boldsymbol{\eta}_{p,m},\boldsymbol{\xi}_{p,m},p_{p,m},\theta_{p,m})$, we can pass to the limit in \eqref{eq:19}--\eqref{eq:25} to obtain a solution to \eqref{eq:28}--\eqref{eq:33}. Reflexivity of $\mathbf{V}_f$, $ T_f$, $\mathbf{H}_p$, $ X_p$ and $T_p$. Consequently, we can extract a subsequence (still denoted by the same symbols) such that

\begin{equation}
\begin{aligned}
\mathbf{u}_{f,m}&\rightharpoonup\mathbf{u}_f &&\text{in }H^1(0,T;\mathbf{H}^1(\Omega_f)),\\
\theta_{f,m}&\rightharpoonup\theta_f &&\text{in }H^1(0,T;H^1(\Omega_f)),\\
\boldsymbol{\xi}_{p,m}&\rightharpoonup\boldsymbol{\xi}_p &&\text{in }H^1(0,T;\mathbf{H}^1(\Omega_p)),\\
\boldsymbol{\eta}_{p,m}&\rightharpoonup\boldsymbol{\eta}_p &&\text{in }H^1(0,T;\mathbf{H}^1(\Omega_p)),\\
p_{p,m}&\rightharpoonup p_p &&\text{in }H^1(0,T;H^1(\Omega_p)),\\
\theta_{p,m}&\rightharpoonup\theta_p &&\text{in }H^1(0,T;H^1(\Omega_p)).
\label{eq:86}
\end{aligned}
\end{equation}
where $\rightharpoonup$ stands for weak convergence. Due to the continuity of the trace operator we also have
By the trace theorem, the traces on $\partial\Omega_{f}$ or $\partial\Omega_{p}$ converge weakly in the corresponding spaces,
\begin{align*}
\mathbf{u}_{f,m}|_{\partial\Omega_{f}}\rightharpoonup\mathbf{u}_f|_{\partial\Omega_{f}}\quad\text{in }L^2(0,T;H^{1/2}(\partial\Omega_{f})),
\end{align*}
and similarly for $\theta_{f,m}$, $\boldsymbol{\xi}_{p,m}$, $p_{p,m}$, $\theta_{p,m}$.

To handle the nonlinear convective terms, we need stronger convergence. Apply \cref{lem:compactness}. For the fluid velocity, take $X=\mathbf{H}^1(\Omega_f)$, $B=\mathbf{L}^4(\Omega_f)$ (since $H^1\hookrightarrow L^4$ is compact in $\mathbb{R}^d$, $d=2,3$), and $Y=\mathbf{L}^2(\Omega_f)$. The sequence $\{\mathbf{u}_{f,m}\}$ is bounded in $L^2(0,T;X)$ and $\{\partial_t\mathbf{u}_{f,m}\}$ is bounded in $L^2(0,T;Y)$. Hence a subsequence (still denoted by $\mathbf{u}_{f,m}$) satisfies
\begin{align}
\mathbf{u}_{f,m}\to\mathbf{u}_f\quad\text{strongly in }L^2(0,T;\mathbf{L}^4(\Omega_f)).\label{eq:87}
\end{align}
Similarly, for the fluid temperature, with $X=H^1(\Omega_f)$, $B=L^4(\Omega_f)$, $Y=L^2(\Omega_f)$, we obtain
\begin{align}
\theta_{f,m}\to\theta_f\quad\text{strongly in }L^2(0,T;L^4(\Omega_f)).\label{eq:88}
\end{align}

All linear terms converge by the weak convergences in \eqref{eq:86}. For the nonlinear term in the fluid momentum equation, we have
\begin{align*}
\int_0^T\big((\mathbf{u}_{f,m}\cdot\nabla)\mathbf{u}_{f,m},\psi\mathbf{v}\big)_{\Omega_f}dt
\to\int_0^T\big((\mathbf{u}_f\cdot\nabla)\mathbf{u}_f,\psi\mathbf{v}\big)_{\Omega_f}dt,
\end{align*}
because $\mathbf{u}_{f,m}\to\mathbf{u}_f$ in $L^2(0,T;\mathbf{L}^4(\Omega_f))$ and $\nabla\mathbf{u}_{f,m}\rightharpoonup\nabla\mathbf{u}_f$ in $L^2(0,T;\mathbf{L}^2(\Omega_f))$ imply $\mathbf{u}_{f,m}\cdot\nabla\mathbf{u}_{f,m}\rightharpoonup\mathbf{u}_f\cdot\nabla\mathbf{u}_f$ in $L^1(0,T;\mathbf{L}^{4/3}(\Omega_f))$. Similarly, by \eqref{eq:87} and \eqref{eq:88},
\begin{align*}
\int_0^T\big((\mathbf{u}_{f,m}\cdot\nabla)\theta_{f,m},\psi w\big)_{\Omega_f}dt
\to\int_0^T\big((\mathbf{u}_f\cdot\nabla)\theta_f,\psi w\big)_{\Omega_f}dt.
\end{align*}

Fix an integer $k\ge1$ and choose $\psi\in L^2(0,T)$. Multiply the Galerkin equations \eqref{eq:28}--\eqref{eq:33} by $\psi(t)$ and integrate over $(0,T)$. For $m\ge k$, using the convergences \eqref{eq:86}--\eqref{eq:88}, we pass to the limit $m\to\infty$.
\begin{align}
&\int_0^T (\partial_t\mathbf{u}_f,\psi(t)\mathbf{v})_{\Omega_f} dt
+\int_0^T 2\mu_f(\mathbb{D}(\mathbf{u}_f),\psi(t)\mathbb{D}(\mathbf{v}))_{\Omega_f} dt
+\int_0^T \langle p_p,\psi(t)\mathbf{v}\cdot\mathbf{n}_f\rangle_\Gamma dt \nonumber\\
&\quad+\int_0^T \frac{\alpha\sqrt{d}}{\sqrt{\operatorname{tr}(\mathbf{K})}}\langle P_\tau(\mathbf{u}_f-\boldsymbol{\xi}_p),P_\tau(\psi(t)\mathbf{v})\rangle_\Gamma dt
+\int_0^T ((\mathbf{u}_f\cdot\nabla)\mathbf{u}_f,\psi(t)\mathbf{v})_{\Omega_f} dt \nonumber\\
&\quad+\int_0^T (\partial_t\theta_f,\psi(t)w)_{\Omega_f} dt +\int_0^T a_f(\nabla\theta_f,\psi(t)\nabla w)_{\Omega_f} dt
+\int_0^T ((\mathbf{u}_f\cdot\nabla)\theta_f,\psi(t)w)_{\Omega_f} dt \nonumber\\
&\quad-\int_0^T a_f\langle\nabla\theta_f\cdot\mathbf{n}_f,\psi(t)w\rangle_\Gamma dt +\int_0^T (\partial_t\boldsymbol{\xi}_p,\psi(t)\boldsymbol{\varphi}_p)_{\Omega_p} dt
+\int_0^T 2\mu_p(\mathbb{D}(\boldsymbol{\eta}_p),\psi(t)\mathbb{D}(\boldsymbol{\varphi}_p))_{\Omega_p} dt \nonumber\\
&\quad+\int_0^T \lambda_p(\nabla\cdot\boldsymbol{\eta}_p,\psi(t)\nabla\cdot\boldsymbol{\varphi}_p)_{\Omega_p} dt -\int_0^T \alpha(\nabla\cdot\boldsymbol{\varphi}_p,\psi(t)p_p)_{\Omega_p} dt
-\int_0^T \beta(\nabla\cdot\boldsymbol{\varphi}_p,\psi(t)\theta_p)_{\Omega_p} dt \nonumber\\
&\quad-\int_0^T \langle p_p,\psi(t)\phi_p\cdot\mathbf{n}_f\rangle_\Gamma dt -\int_0^T \frac{\alpha\sqrt{d}}{\sqrt{\operatorname{tr}(\mathbf{K})}}\langle P_\tau(\mathbf{u}_f-\boldsymbol{\xi}_p),P_\tau(\psi(t)\phi_p)\rangle_\Gamma dt \nonumber\\
&\quad+\int_0^T S_0(\partial_t p_p,\psi(t)q_p)_{\Omega_p} dt
-\int_0^T b_0(\partial_t\theta_p,\psi(t)q_p)_{\Omega_p} dt +\int_0^T \alpha(\nabla\cdot\boldsymbol{\xi}_p,\psi(t)q_p)_{\Omega_p} dt \nonumber\\
&\quad+\int_0^T \mathbf{K}(\nabla p_p,\psi(t)\nabla q_p)_{\Omega_p} dt
-\int_0^T \langle(\mathbf{u}_f-\boldsymbol{\xi}_p)\cdot\mathbf{n}_f,\psi(t)q_p\rangle_\Gamma dt \nonumber\\
&\quad+\int_0^T a_0(\partial_t\theta_p,\psi(t)w_p)_{\Omega_p} dt
-\int_0^T b_0(\partial_t p_p,\psi(t)w_p)_{\Omega_p} dt
+\int_0^T \beta(\nabla\cdot\boldsymbol{\xi}_p,\psi(t)w_p)_{\Omega_p} dt \nonumber\\
&\quad+\int_0^T \bm{a}_p(\nabla\theta_p,\psi(t)\nabla w_p)_{\Omega_p} dt
-\int_0^T \bm{a}_p\langle\nabla\theta_p\cdot\mathbf{n}_p,\psi(t)w_p\rangle_\Gamma dt \nonumber\\
&= \int_0^T \frac{Gr}{Re^2}(\boldsymbol{\xi}\theta_f,\psi(t)\mathbf{v})_{\Omega_f} dt
+\int_0^T (\mathbf{u}_f\cdot\boldsymbol{\xi},\psi(t)w)_{\Omega_f} dt
+\int_0^T (\mathbf{f},\psi(t)\phi_p)_{\Omega_p} dt \nonumber\\
&\quad+\int_0^T (g,\psi(t)q_p)_{\Omega_p} dt
+\int_0^T (\phi,\psi(t)w_p)_{\Omega_p} dt,\label{eq:89}
\end{align}

for all $(\mathbf{v},w,\bm{\varphi_p},q_p,w_p)\in\mathbf{V}_f^k\times T_f^k\times\mathbf{H}_p^k\times X_p^k\times T_p^k$. Since $k$ is arbitrary and the union of these finite-dimensional spaces is dense in the respective spaces, and $\psi(t)\in L^2(0,T)$ is arbitrary, we conclude that the equality holds for almost every $t\in(0,T)$ and for all test functions in $\mathbf{V}_f\times T_f\times\mathbf{H}_p\times X_p\times T_p$. Thus $(\mathbf{u}_f,\theta_f,\boldsymbol{\eta}_p,\boldsymbol{\xi}_p,p_p,\theta_p)$ satisfies \eqref{eq:28}--\eqref{eq:33}.

It remains to verify the initial conditions. Choose $\psi(t)\in C^2([0,T])$ with $\psi(T)=\partial_t \psi(T)=0$ and integrate by parts the time derivative terms .
\begin{align*}
\int_0^T(\partial_t\mathbf{u}_{f,m},\psi(t)\mathbf{v})dt
=-\int_0^T(\mathbf{u}_{f,m},\partial_t \psi(t)\mathbf{v})dt-(\mathbf{u}_{f,m}(0),\psi(0)\mathbf{v}).
\end{align*}
Proceed similarly for $\partial_t\theta_{f,m}$, $\partial_t\boldsymbol{\xi}_{p,m}$, $\partial_t p_{p,m}$, $\partial_t\theta_{p,m}$. Letting $m\to\infty$ and using the convergences \eqref{eq:86}--\eqref{eq:88}, we obtain
\begin{align*}
&-\int_0^T(\mathbf{u}_f,\partial_t \psi(t)\mathbf{v})dt-(\mathbf{u}_f(0),\psi(0)\mathbf{v})-\int_0^T(\theta_f,\partial_t \psi(t)w)dt-(\theta_f(0),\psi(0)w)\\
&-\int_0^T(\boldsymbol{\xi}_p,\partial_t \psi(t)\phi_p)dt-(\boldsymbol{\xi}_p(0),\psi(0)\phi_p)+\int_0^T(\boldsymbol{\eta}_p,\partial_t^2{\psi}(t)\phi_p)dt+(\boldsymbol{\eta}_p(0),\partial_t \psi(t)(0)\phi_p)\\
&-\int_0^T(S_0 p_p,\partial_t \psi(t)q_p)dt-(S_0 p_p(0),\psi(0)q_p)+\int_0^T(b_0\theta_p,\partial_t \psi(t)q_p)dt+(b_0\theta_p(0),\psi(0)q_p)\\
&-\int_0^T(a_0\theta_p,\partial_t \psi(t)w_p)dt-(a_0\theta_p(0),\psi(0)w_p)+\int_0^T(b_0 p_p,\partial_t \psi(t)w_p)dt+(b_0 p_p(0),\psi(0)w_p)= 0,
\end{align*}
then
\begin{align*}
&-(\mathbf{u}_f(0),\psi(0)\mathbf{v})-(\theta_f(0),\psi(0)w)-(\boldsymbol{\xi}_p(0),\psi(0)\phi_p)+(\boldsymbol{\eta}_p(0),\partial_t \psi(0)\phi_p)\\
&-(S_0 p_p(0),\psi(0)q_p)+(b_0\theta_p(0),\psi(0)q_p)-(a_0\theta_p(0),\psi(0)w_p)+(b_0 p_p(0),\psi(0)w_p)= 0,
\end{align*}
for all $(\mathbf{v},w,\bm{\varphi}_p,q_p,w_p)\in\mathbf{V}_f^k\times T_f^k\times\mathbf{H}_p^k\times X_p^k\times T_p^k$. and for all $\psi(0),\partial_t \psi(0)$  are arbitrary, we deduce
\begin{align*}
\mathbf{u}_f(0)=0,\;\theta_f(0)=0,\;\boldsymbol{\xi}_p(0)=0,\;\boldsymbol{\eta}_p(0)=0,\;\;p_p(0)=0,\;\theta_p(0)=0.
\end{align*}
This yields the initial conditions stated in \eqref{eq:18}. Finally, passing to the limit in \eqref{eq:42}--\eqref{eq:44}, we complete the proof of existence of a weak solution $(\mathbf{u}_{f},\theta_{f},\boldsymbol{\eta}_{p},\boldsymbol{\xi}_{p},p_{p},\theta_{p})$ to \eqref{eq:12}--\eqref{eq:18}.

\section{Uniqueness analysis}\label{sec:4}
Having established the existence of a weak solution, we next turn to the question of uniqueness. The following theorem shows that, under a smallness condition analogous to that used for existence, the solution is indeed unique, thereby completing the well-posedness analysis.

\begin{theorem}
\label{thm:local_uniqueness}
Let $(\mathbf{u}_{f},\theta_{f},\boldsymbol{\eta}_{p},\boldsymbol{\xi}_{p},p_{p},\theta_{p})$ be a solution of \eqref{eq:12}--\eqref{eq:18}. Then if
\begin{align}
\|\mathbb{D}(\mathbf{u}_f)\|_{L^2} \le \min\left(\frac{3\mu_f}{4C_{S,f}^2 C_{K}^3},\ \frac{a_f}{2C_{S_f}^2 C_K C_{P,f} }\right),\|\nabla\theta_{f}\|_{L^2} \le \frac{\sqrt{a_f \mu_f}}{C_{S_f}^2 C_K C_{P,f} }
\label{eq:900}
\end{align}
then $(\mathbf{u}_{f},\theta_{f},\boldsymbol{\eta}_{p},\boldsymbol{\xi}_{p},p_{p},\theta_{p})$ is unique.
\end{theorem}
\begin{proof}\
Let $(\mathbf{u}_{f,1},\theta_{f,1},\boldsymbol{\eta}_{p,1},\boldsymbol{\xi}_{p,1},p_{p,1},\theta_{p,1})$ and $(\mathbf{u}_{f,2},\theta_{f,2},\boldsymbol{\eta}_{p,2},\boldsymbol{\xi}_{p,2},p_{p,2},\theta_{p,2})$ be two solutions of \eqref{eq:12}--\eqref{eq:18} with the same initial and boundary data.
Define the differences
\begin{align*}
e_{\mathbf{u}} &= \mathbf{u}_{f,1}-\mathbf{u}_{f,2}, \quad
e_{\theta_f} = \theta_{f,1}-\theta_{f,2}, \quad
e_{\bm{\eta}}= \bm{\eta}_{p,1}-\bm{\eta}_{p,2}, \\
e_{\bm{\xi}} &= \bm{\xi}_{p,1}-\bm{\xi}_{p,2} = \partial_t \mathbf{e}_{\bm{\eta}}, \quad
e_p = p_{p,1}-p_{p,2}, \quad
e_{\theta p} = \theta_{p,1}-\theta_{p,2}.
\end{align*}
The initial conditions imply
\begin{align*}
e_u(0)=\mathbf{0},\
e_{\theta_f}(0)=0,\
e_{\eta}(0)=\mathbf{0},\
e_{\xi}(0)=\mathbf{0},\
e_p(0)=0,\
e_{\theta _p}(0)=0.
\end{align*}
Subtracting the weak formulations for the two solutions yields the following five equations :
\begin{align}
&(\partial_te_u,\mathbf{v})_{\Omega_f}+2\mu_f(\mathbb{D}(e_u),\mathbb{D}(\mathbf{v}))_{\Omega_f}+\big((\mathbf{u}_{f,1}\cdot\nabla)\mathbf{u}_{f,1}-(\mathbf{u}_{f,2}\cdot\nabla)\mathbf{u}_{f,2},\mathbf{v}\big)_{\Omega_f} \nonumber \\
&\quad+\langle e_p,\mathbf{v}\cdot\mathbf{n}_f\rangle_\Gamma+\frac{\alpha\sqrt{d}}{\sqrt{\operatorname{tr}(\mathbf{K})}}\langle P_\tau(e_u-e_\xi),P_\tau(\mathbf{v})\rangle_\Gamma=0 \label{eq:90},\\[4pt]
&(\partial_t e_{\theta_ f},w)_{\Omega_f}+a_f(\nabla e_{\theta _f},\nabla w)_{\Omega_f}+\big((\mathbf{u}_{f,1}\cdot\nabla)\theta_{f,1}-(\mathbf{u}_{f,2}\cdot\nabla)\theta_{f,2},w\big)_{\Omega_f} \nonumber \\
&\quad-a_f\langle\nabla e_{\theta f}\cdot\mathbf{n}_f,w\rangle_\Gamma=0 \label{eq:91},\\[4pt]
&(\partial_te_\xi,\boldsymbol{\varphi}_p)_{\Omega_p}+2\mu_p(\mathbb{D}(e_\eta),\mathbb{D}(\boldsymbol{\varphi}_p))_{\Omega_p}+\lambda_p(\nabla\cdot e_\eta,\nabla\cdot\boldsymbol{\varphi}_p)_{\Omega_p}-\alpha(\nabla\cdot\boldsymbol{\varphi}_p,e_p)_{\Omega_p} \nonumber \\
&\quad-\beta(\nabla\cdot\boldsymbol{\varphi}_p,e_{\theta p})_{\Omega_p}-\langle e_p,\boldsymbol{\varphi}_p\cdot\mathbf{n}_f\rangle_\Gamma-\frac{\alpha\sqrt{d}}{\sqrt{\operatorname{tr}(\mathbf{K})}}\langle P_\tau(e_u-e_\xi),P_\tau(\boldsymbol{\varphi}_p)\rangle_\Gamma=0\label{eq:92}, \\[4pt]
&S_0(\partial_t e_p,q_p)_{\Omega_p}-b_0(\partial_t e_{\theta_ p},q_p)_{\Omega_p}+\alpha(\nabla\cdot e_\xi,q_p)_{\Omega_p}+\mathbf{K}(\nabla e_p,\nabla q_p)_{\Omega_p} \nonumber \\
&\quad-\langle(e_u-e_\xi)\cdot\mathbf{n}_f,q_p\rangle_\Gamma=0 \label{eq:93},\\[4pt]
&a_0(\partial_t e_{\theta_ p},w_p)_{\Omega_p}-b_0(\partial_t e_p,w_p)_{\Omega_p}+\beta(\nabla\cdot\mathbf{e}_\xi,w_p)_{\Omega_p}+\bm{a}_p(\nabla e_{\theta_p},\nabla w_p)_{\Omega_p} \nonumber \\
&\quad-\bm{a}_p\langle\nabla e_{\theta p}\cdot\mathbf{n}_p,w_p\rangle_\Gamma=0\label{eq:94}.
\end{align}
Now choose the test functions as
$\mathbf{v}=e_u$, $w=e_{\theta_f}$, $\bm{\varphi}_p=\mathbf{e}_\xi$, $q_p=e_p$, $w_p=e_{\theta_p}$. Insert these into \eqref{eq:90}--\eqref{eq:94} and add the resulting equations,
 we obtain
\begin{equation}
\begin{aligned}
&\frac{1}{2}\frac{d}{dt}\|e_u\|^2_{L^2}+\frac{1}{2}\frac{d}{dt}\|e_{\theta_ f}\|^2_{L^2}+\frac{1}{2}\frac{d}{dt}\|e_\xi\|^2_{L^2}
+\mu_p\frac{d}{dt}\|\mathbb{D}(e_\eta)\|^2_{L^2}+\frac{\lambda_p}{2}\frac{d}{dt}\|\nabla\cdot e_\eta\|^2_{L^2}\\
&+\frac{S_0}{2}\frac{d}{dt}\|e_p\|^2_{L^2}+\frac{a_0}{2}\frac{d}{dt}\|e_{\theta_ p}\|^2_{L^2}-b_0\frac{d}{dt}(e_p,e_{\theta_ p})+2\mu_f\|\mathbb{D}(e_u)\|^2_{L^2}+a_f\|\nabla e_{\theta_ f}\|^2_{L^2}\\
&+\mathbf{K}\|\nabla e_p\|^2_{L^2}+\bm{a}_p\|\nabla e_{\theta p}\|^2_{L^2}+\frac{\alpha\sqrt{d}}{\sqrt{\operatorname{tr}(\mathbf{K})}} \| P_\tau(e_u - e_\xi) \|_{\Gamma}^2\\
&=-I_u-I_\theta,
\end{aligned}
\label{eq:energy}
\end{equation}
where
\[
I_u:=\big((\mathbf{u}_{f,1}\cdot\nabla)\mathbf{u}_{f,1}-(\mathbf{u}_{f,2}\cdot\nabla)\mathbf{u}_{f,2},\mathbf{e}_u\big)_{\Omega_f},\quad
I_\theta:=\big((\mathbf{u}_{f,1}\cdot\nabla)\theta_{f,1}-(\mathbf{u}_{f,2}\cdot\nabla)\theta_{f,2},e_{\theta_f}\big)_{\Omega_f}.
\]
Now we estimate $I_u$ and $I_\theta$. Using the decompositions
\[
(\mathbf{u}_{f,1}\cdot\nabla)\mathbf{u}_{f,1}-(\mathbf{u}_{f,2}\cdot\nabla)\mathbf{u}_{f,2}=(e_u\cdot\nabla)\mathbf{u}_{f,1}+(\mathbf{u}_{f,2}\cdot\nabla)e_u,
\]
\[
(\mathbf{u}_{f,1}\cdot\nabla)\theta_{f,1}-(\mathbf{u}_{f,2}\cdot\nabla)\theta_{f,2}=(e_u\cdot\nabla)\theta_{f,1}+(\mathbf{u}_{f,2}\cdot\nabla)e_{\theta_f},
\]
and applying Hölder's inequality, the Sobolev embedding $H^1\hookrightarrow L^4$ , Korn's inequality and Young's inequality, we obtain
\begin{equation}
\begin{aligned}
|I_u| &\le C_{S,f}^2 C_{K}^3 \|\mathbb{D}(e_u)\|_{L^2}^2 \|\mathbb{D}(\mathbf{u}_{f,1})\|_{L^2} + C_{S,f}^2 C_{K}^3 \|\mathbb{D}(e_u)\|_{L^2}^2 \|\mathbb{D}(\mathbf{u}_{f,2})\|_{L^2}, \\
|I_\theta| &\le C_{S_f}^2 C_K C_{p,f} \|\nabla\theta_{f,1}\|_{L^2} \|\mathbb{D}(e_u)\|_{L^2} \|\nabla e_{\theta f}\|_{L^2} + C_{S_f}^2 C_K C_{p,f} \|\mathbb{D}(\mathbf{u}_{f,2})\| _{L^2}\|\nabla e_{\theta_ f}\|^2_{L^2}\\
&
\le\|\nabla\theta_{f,1}\|_{L^2} \frac{\mu_f}{2} \|\mathbb{D}(e_u)\|^2_{L^2} + \|\nabla\theta_{f,1}\|_{L^2}\frac{ (C_{S_f}^2 C_K C_{p,f})^2 }{2\mu_f} \|\nabla e_{\theta_ f}\|^2_{L^2}\\&+ C_{S_f}^2 C_K C_{p,f} \|\mathbb{D}(\mathbf{u}_{f,2})\|_{L^2} \|\nabla e_{\theta_ f}\|^2_{L^2}.
\end{aligned}
\label{eq:estimates}
\end{equation}
Consequently, using \eqref{eq:900}, we can obtain the following estimate:
\begin{align}
|I_u|+|I_\theta| \le 2\mu_f \|\mathbb{D}(\mathbf{e}_u)\|^2_{L^2} + a_f \|\nabla e_{\theta_ f}\|^2_{L^2}.\label{eq:98}
\end{align}

Due to the thermo‑poroelastic stability condition $S_0a_0 > b_0^2$, there exists a constant $c_0>0$ such that
\begin{align}
\frac{S_0}{2}\|e_p\|^2+\frac{a_0}{2}\|e_{\theta p}\|^2-b_0(e_p,e_{\theta p})\ge\frac{c_0}{2}\bigl(\|e_p\|^2+\|e_{\theta p}\|^2\bigr),
\label{eq:99}
\end{align}
hence the time derivative part can be bounded below by $\frac{c_0}{2}\frac{d}{dt}(\|e_p\|^2+\|e_{\theta p}\|^2)$.

Substituting \eqref{eq:98} and \eqref{eq:99} into \eqref{eq:energy} yields
\begin{align*}
&\frac12\frac{d}{dt}\|e_u\|^2_{L^2} + \frac12\frac{d}{dt}\|e_{\theta_ f}\|^2_{L^2} + \frac12\frac{d}{dt}\|e_\xi\|^2_{L^2}
+ \frac{c_0}{2}\frac{d}{dt}\|e_p\|^2_{L^2} + \frac{c_0}{2}\frac{d}{dt}\|e_{\theta_ p}\|^2_{L^2} \\
&\quad + \mu_p\frac{d}{dt}\|\mathbb{D}(e_\eta)\|^2_{L^2} + \frac{\lambda_p}{2}\frac{d}{dt}\|\nabla\cdot e_\eta\|^2_{L^2}
\le 0.
\end{align*}
By the nonnegativity of all terms on the left-hand side and the initial conditions, it follows that
\begin{align*}
\frac12\|e_u\|^2_{L^2} + \frac12\|e_{\theta_ f}\|^2_{L^2} + \frac12\|e_\xi\|^2_{L^2}
+ \frac{c_0}{2}\|e_p\|^2_{L^2} + \frac{c_0}{2}\|e_{\theta_ p}\|^2_{L^2}
+ \mu_p\|\mathbb{D}(e_\eta)\|^2_{L^2} + \frac{\lambda_p}{2}\|\nabla\cdot e_\eta\|^2_{L^2}=0.
\end{align*} Consequently,
\begin{align*}
e_u \equiv \bm0,\quad e_{\theta_ f} \equiv 0,\quad e_\xi \equiv \bm0,\quad e_p \equiv 0,\quad e_{\theta_ p} \equiv 0,\quad e_\eta \equiv \bm0,
\end{align*}
which implies $\mathbf{u}_{f,1}=\mathbf{u}_{f,2}$, $\theta_{f,1}=\theta_{f,2}$, $\boldsymbol{\eta}_{p,1}=\boldsymbol{\eta}_{p,2}$, $p_{p,1}=p_{p,2}$, $\theta_{p,1}=\theta_{p,2}$.
\end{proof}
% ============================================================

\begin{remark}\label{rem:data-control}
Condition \eqref{eq:900} may appear to be a control on the
solution $(\mathbf u_f,\theta_f)$ rather than on the data. In fact, it can be
regarded as an equivalent way of writing a smallness assumption on the data.
The uniform estimates established in \cref{sec:3} for the Galerkinc
approximations hold with constants independent of $m$; passing to the limit
$m\to\infty$ as in \cref{subsec:limit}  yields the same type of estimates for the limit
solution. Consequently,
\eqref{eq:900} follows from the purely data-dependent
smallness of $\mathcal C_1,\mathcal C_2,\mathcal C_3$ defined in
\eqref{eq:a30}, and therefore does not constitute an additional restriction
on the solution.
\end{remark}
\section{Recovery of the fluid pressure}\label{sec:5}

The solution constructed so far satisfies a weak formulation in which the fluid pressure
$p_f$
  has been eliminated by restricting to divergence-free test functions. In this final step, we recover the pressure via the classical inf--sup condition, thereby showing that the solution also satisfies the original coupled system \eqref{eq:12}--\eqref{eq:18} with a uniquely defined fluid pressure.

\begin{theorem}
\label{thm:existence_pf}
If $(\mathbf{u}_f,\theta_f,\boldsymbol{\eta}_p,\boldsymbol{\xi}_p,p_p,\theta_p)$ is a solution of problem \eqref{eq:19}--\eqref{eq:25}, then there exists a unique $p_f \in L^\infty(0,T; X_f)$ such that $(\mathbf{u}_f, p_f, \theta_f,\boldsymbol{\eta}_p,\boldsymbol{\xi}_p,p_p,\theta_p)$ is a solution of problem \eqref{eq:12}--\eqref{eq:18}.
\end{theorem}
\begin{proof}\
Define a linear functional $\mathcal{F}$ on $\mathbf H_f$ by
\begin{align*}
&\mathcal{F}(\mathbf v_f)
= (\partial_t\mathbf u_f,\mathbf v_f)_{\Omega_f}
+ 2\mu_f(\mathbb D(\mathbf u_f),\mathbb D(\mathbf v_f))_{\Omega_f}
+ ((\mathbf u_f\cdot\nabla)\mathbf u_f,\mathbf v_f)_{\Omega_f} \notag\\
&\quad - \frac{Gr}{Re^2}(\bm{\xi}\theta_f,\mathbf v_f)_{\Omega_f}
+ \langle p_p,\mathbf v_f\cdot\mathbf n_f\rangle_\Gamma
+ \frac{\alpha\sqrt d}{\sqrt{\operatorname{tr}(\mathbf K)}}
\langle P_\tau(\mathbf u_f-\bm{\xi}_p),P_\tau(\mathbf v_f)\rangle_\Gamma.
\end{align*}

 For any $\mathbf{v}_f\in\mathbf{V}_f$ (the subspace of divergence-free functions), the equations \eqref{eq:19}--\eqref{eq:25} imply $\mathcal{F}(\mathbf{v}_f)=0$.

The fluid domain $\Omega_f$ satisfies the inf-sup condition \cite{cesmelioglu2017analysis}: there exists a constant $L>0$ such that
\begin{align}
\inf_{q\in X_f}\sup_{\mathbf{v}_f\in\mathbf{H}_f}
\frac{(\nabla\cdot\mathbf{v}_f,q)_{\Omega_f}}{\|\mathbf{v}_f\|_{H^1(\Omega_f)}\|q\|_{L^2(\Omega_f)}}\ge L.
\end{align}
Since $\mathcal{F}$ is a continuous linear functional that vanishes on $\mathbf H_f$, the inf-sup condition guarantees the existence of a unique $p_f\in X_f$ for almost every $t\in(0,T)$ such that
\begin{align}
\mathcal{F}(\mathbf{v}_f)=(\nabla\cdot\mathbf{v}_f,p_f)_{\Omega_f},\qquad\forall\,\mathbf{v}_f\in\mathbf{H}_f,
\end{align}
and we have the estimate
\begin{align}
\|p_f\|_{L^2(\Omega_f)}\le\frac{1}{L}\sup_{\mathbf{v}_f\in\mathbf{H}_f}
\frac{|\mathcal{F}(\mathbf{v}_f)|}{\|\mathbf{v}_f\|_{H^1(\Omega_f)}}.\label{eq:102}
\end{align}

We bound the right-hand side of \eqref{eq:102} using the a priori estimates already established for the solution of \eqref{eq:19}--\eqref{eq:25}. All constants below are independent of the solution; they include the Poincaré constant $C_{P,f}$, the Korn constant $C_{K,f}$, the Sobolev embedding constant $C_{S,f}$ and the trace constants $T_1,T_3,T_5$ introduced in \cref{sec:6}.

\begin{gather}
|(\partial_t\mathbf{u}_f,\mathbf{v}_f)|
\le \|\partial_t\mathbf{u}_f\|_{L^2}\|\mathbf{v}_f\|_{L^2}
\le C_{P,f}\|\partial_t\mathbf{u}_f\|_{L^2}\|\nabla\mathbf{v}_f\|_{L^2} \le C_{P,f}C_{K,f}\|\partial_t\mathbf{u}_f\|_{L^2}\|\mathbb{D}(\mathbf{v}_f)\|_{L^2},
\label{eq:103}
\\[4pt]
|2\mu_f(\mathbb{D}(\mathbf{u}_f),\mathbb{D}(\mathbf{v}_f))|
\le 2\mu_f\|\mathbb{D}(\mathbf{u}_f)\|_{L^2}\|\mathbb{D}(\mathbf{v}_f)\|_{L^2},\label{eq:104}
\\[4pt]
|((\mathbf{u}_f\cdot\nabla)\mathbf{u}_f,\mathbf{v}_f)|
\le \|\mathbf{u}_f\|_{L^4}\|\nabla\mathbf{u}_f\|_{L^2}\|\mathbf{v}_f\|_{L^4} \le C_{S,f}^2\|\mathbf{u}_f\|_{H^1}^2\|\mathbf{v}_f\|_{H^1}
\le C_{S,f}^2C_{K,f}^2\|\mathbb{D}(\mathbf{u}_f)\|_{L^2}^2\|\mathbf{v}_f\|_{H^1},
\label{eq:105}\\[4pt]
|\langle p_p,\mathbf{v}_f\cdot\mathbf{n}_f\rangle_\Gamma|
\le \|p_p\|_{L^2(\Gamma)}\|\mathbf{v}_f\|_{L^2(\Gamma)}
\le T_3\|p_p\|_{H^1(\Omega_p)}\cdot T_1\|\mathbf{v}_f\|_{H^1(\Omega_f)},\label{eq:106}
\end{gather}
\begin{equation}
\begin{split}
&|\frac{\alpha\sqrt{d}}{\sqrt{\operatorname{tr}(\mathbf{K})}}\langle P_\tau(\mathbf{u}_f-\boldsymbol{\xi}_p),P_\tau(\mathbf{v}_f)\rangle_\Gamma|
\le \frac{\alpha\sqrt{d}}{\sqrt{\operatorname{tr}(\mathbf{K})}}\bigl(\|\mathbf{u}_f\|_{L^2(\Gamma)}+\|\boldsymbol{\xi}_p\|_{L^2(\Gamma)}\bigr)\|\mathbf{v}_f\|_{L^2(\Gamma)} \\
&\le \frac{\alpha\sqrt{d}}{\sqrt{\operatorname{tr}(\mathbf{K})}}\bigl(T_1\|\mathbf{u}_f\|_{H^1}+T_5\|\boldsymbol{\xi}_p\|_{H^1}\bigr)\,T_1\|\mathbf{v}_f\|_{H^1} \\
&\le \frac{\alpha\sqrt{d}}{\sqrt{\operatorname{tr}(\mathbf{K})}} T_1^2C_{K,f}\|\mathbb{D}(\mathbf{u}_f)\|_{L^2}\|\mathbf{v}_f\|_{H^1}
   + \frac{\alpha\sqrt{d}}{\sqrt{\operatorname{tr}(\mathbf{K})}}  T_1T_5\|\boldsymbol{\xi}_p\|_{H^1}\|\mathbf{v}_f\|_{H^1},\label{eq:107}
\end{split}
\end{equation}
\begin{align}
\Bigl|\frac{G_r}{Re^2}(\boldsymbol{\xi}\theta_f,\mathbf{v}_f)\Bigr|
\le \frac{G_r}{Re^2}\|\theta_f\|_{L^2}\|\mathbf{v}_f\|_{L^2}
\le \frac{G_r}{Re^2}C_{P,f}\|\theta_f\|_{L^2}\|\mathbf{v}_f\|_{H^1}.\label{eq:108}
\end{align}

Collecting the estimates \eqref{eq:103}--\eqref{eq:108}and noting that $\|\mathbb{D}(\mathbf{v}_f)\|_{L^2}\le\|\mathbf{v}_f\|_{H^1}$, we obtain

\begin{equation}
\begin{split}
\frac{|\mathcal{F}(\mathbf{v}_f)|}{\|\mathbf{v}_f\|_{H^1}}
&\le C_{P,f}C_{K,f}\|\partial_t\mathbf{u}_f\|_{L^2}
+ 2\mu_f\|\mathbb{D}(\mathbf{u}_f)\|_{L^2}
+ C_{S,f}^2 C_{K,f}^2 \|\mathbb{D}(\mathbf{u}_f)\|_{L^2}^2  + T_1 T_3 \|p_p\|_{H^1}
\\
&\quad+ \frac{\alpha\sqrt{d}}{\sqrt{\operatorname{tr}(\mathbf{K})}}  T_1^2 C_{K,f} \|\mathbb{D}(\mathbf{u}_f)\|_{L^2}  + \frac{\alpha\sqrt{d}}{\sqrt{\operatorname{tr}(\mathbf{K})}}  T_1 T_5 \|\boldsymbol{\xi}_p\|_{H^1}
+ \frac{G_r}{Re^2}C_{P,f} \|\theta_f\|_{L^2}.
\end{split}
\end{equation}
 Consequently,
\begin{equation}
\begin{split}
\|p_f\|_{L^2}
\le \frac{1}{L}\Bigl( & C_{P,f}C_{K,f}\|\partial_t\mathbf{u}_f\|_{L^2}
+ 2\mu_f\|\mathbb{D}(\mathbf{u}_f)\|_{L^2}
+ C_{S,f}^2C_{K,f}^2\|\mathbb{D}(\mathbf{u}_f)\|_{L^2}^2 + T_1T_3\|p_p\|_{H^1}
\\
& + \frac{\alpha\sqrt{d}}{\sqrt{\operatorname{tr}(\mathbf{K})}} T_1^2C_{K,f}\|\mathbb{D}(\mathbf{u}_f)\|_{L^2} + \frac{\alpha\sqrt{d}}{\sqrt{\operatorname{tr}(\mathbf{K})}} T_1T_5\|\boldsymbol{\xi}_p\|_{H^1}
+ \frac{G_r}{Re^2}C_{P,f}\|\theta_f\|_{L^2}\Bigr).
\end{split}
\end{equation}
Furthermore, using the \(m\)-uniform estimates \eqref{eq:42}--\eqref{eq:43}, we get
\begin{equation}
\begin{split}
\|\partial_t\mathbf{u}_f(t)\|_{L^2(\Omega_f)}^2
&\le 2(1+Te^T)\Big(\|\mathcal C_2\|_{L^2(0,T)}^2+\mathcal C_3^2\Big),\\
\|p_p(t)\|_{H^1(\Omega_p)}^2
&\le \frac{2(1+C_{P,p}^2)(1+e^T)}{K_{min}}\|\mathcal C_1\|_{L^2(0,T)}^2,\\
\|\boldsymbol{\xi}_p(t)\|_{H^1(\Omega_p)}^2
&\le \frac{2(1+C_{P,p}^2)C_K^2(1+Te^T)}{\mu_p}\Big(\|\mathcal C_2\|_{L^2(0,T)}^2+\mathcal C_3^2\Big),
\end{split}
\label{eq:pf-bounds}
\end{equation}
together with \eqref{eq:a31}
we obtain, for a.e.\ $t\in(0,T)$,
\begin{equation}
\begin{split}
\|p_f(t)\|_{L^2(\Omega_f)}
\le &\frac{1}{L}\Bigl(
C_{P,f}C_{K,f}\sqrt{2(1+Te^T)\big(\|\mathcal C_2\|_{L^2(0,T)}^2+\mathcal C_3^2\big)}
+ 2\mu_f\,\delta_1
+ C_{S,f}^2C_{K,f}^2\,\delta_1^2
\\
& + T_1T_3\sqrt{\frac{2(1+C_{P,p}^2)(1+e^T)}{K_{min}}\|\mathcal C_1\|_{L^2(0,T)}^2}
+ \frac{\alpha\sqrt{d}}{\sqrt{\operatorname{tr}(\mathbf{K})}} T_1^2C_{K,f}\,\delta_1
\\
& + \frac{\alpha\sqrt{d}}{\sqrt{\operatorname{tr}(\mathbf{K})}} T_1T_5
\sqrt{\frac{2(1+C_{P,p}^2)C_K^2(1+Te^T)}{\mu_p}\big(\|\mathcal C_2\|_{L^2(0,T)}^2+\mathcal C_3^2\big)}
 + \frac{G_r}{Re^2}C_{P,f}^2\,\delta_2\Bigr).
\end{split}
\label{eq:pf-explicit}
\end{equation}
Here \(\delta_1,\delta_2\) are given in \eqref{eq:a31} and \(\mathcal C_1,\mathcal C_2,\mathcal C_3\) in \eqref{eq:a30}; the right-hand side depends only on the data $\mathbf f,g,\phi$, on \(T\).

The proof is complete.
\end{proof}
\section{Final Results}\label{sec:6}
In the preceding sections (Section~\ref{sec:2}--Section~\ref{sec:5}), we have presented the weak formulation of the thermal--fluid--solid coupling problem and established the existence, uniqueness, and pressure recovery results. Based on those analyses, this section states the main theoretical result of the paper, i.e., the existence and uniqueness theorems for the weak solution. The proof relies on suitable smallness assumptions, which are explicitly listed in the theorem below together with the a priori estimates satisfied by the solution. Detailed proofs have been provided in Section~\ref{sec:3}--Section~\ref{sec:5} using the Galerkin method, compactness arguments, and the inf-sup condition.
 \begin{theorem}\label{thm:main}
Assume that $\mathbf{f} \in H^1(0,T; L^2(\Omega_f))$, $g \in H^1(0,T; L^2(\Omega_p))$, $\phi \in H^1(0,T; L^2(\Omega_p))$ and that the following small data condition holds:
\begin{equation}
\begin{aligned}
&
\|\mathcal{C}_1\|_{L^\infty(0,T)}^2+(1 + e^{T}+Te^{T})\|\mathcal{C}_1\|_{L^2(0,T)}^2  +  (1+Te^{T})  \|\mathcal{C}_2\|_{L^2(0,T)}^2 +  (1+Te^{T})\mathcal{C}_3^2<\mu_f\delta_1,\\
&
\|\mathcal{C}_1\|_{L^\infty(0,T)}^2+(1 + e^{T}+Te^{T})\|\mathcal{C}_1\|_{L^2(0,T)}^2  +  (1+Te^{T})  \|\mathcal{C}_2\|_{L^2(0,T)}^2 +  (1+Te^{T})\mathcal{C}_3^2<\frac{a_f}{2}\delta_2.
\end{aligned}
\end{equation}
Then, problem \eqref{eq:12}--\eqref{eq:18} has a unique solution $(\mathbf{u}_f, p_f, \theta_f,\bm{\eta}_p,\bm{\xi}_p,\theta_p, p_p)$ such that
\begin{equation}
\begin{aligned}
&\frac{1}{2}\norm{\mathbf{u}_{f}}_{L^\infty(0,T;L^2)}^2 + \frac{1}{2}\norm{\theta_{f}}_{L^\infty(0,T;L^2)}^2 + \frac{1}{2}\norm{\bm{\xi}_{p}}_{L^\infty(0,T;L^2)}^2 + \frac{c_0}{2}\norm{p_{p}}_{L^\infty(0,T;L^2(\Omega_p))}^2 \\
&\quad+ \frac{c_0}{2}\norm{\theta_{p}}_{L^\infty(0,T;L^2)}^2  + \mu_p\norm{\mathbb{D}(\bm{\eta}_{p})}_{L^\infty(0,T;L^2)}^2 + \frac{\lambda_p}{2}\norm{\nabla\cdot\bm{\eta}_{p}}_{L^\infty(0,T;L^2)}^2 \\
&\quad+ \mu_f \norm{\mathbb{D}(\mathbf{u}_{f})}_{L^\infty(0,T;L^2)}^2
+ \frac{a_f}{2}\norm{\nabla\theta_{f}}_{L^\infty(0,T;L^2)}^2 \\
&\quad + \frac{\mathbf{K}}{2}\norm{\nabla p_{p}}_{L^\infty(0,T;L^2)}^2 + \frac{\bm{a}_p}{2}\norm{\nabla\theta_{p}}_{L^\infty(0,T;L^2)}^2 \leq ( 1+e^{T})\|\mathcal{C}_1\|_{L^2(0,T)}^2.
\end{aligned}
\end{equation}
Furthermore,
\begin{align}
\|\mathbb{D}(\mathbf{u}_{f})(t)\|_{L^\infty(0,T;L^2)} < \delta_1,\qquad
\|\nabla\theta_{f}(t)\|_{L^\infty(0,T;L^2)} < \delta_2  ,
\end{align}
and
\begin{equation}
\begin{split}
&\frac12\|\partial_t\mathbf{u}_{f}\|^2_{L^\infty(0,T;L^2)}
+ \frac12\|\partial_t\theta_{f}\|^2_{L^\infty(0,T;L^2)}
+ \frac12\|\partial_t\boldsymbol{\xi}_{p}\|^2_{L^\infty(0,T;L^2)}
+ \frac{c_0}{2}\|\partial_t p\|^2_{L^\infty(0,T;L^2)}
\\
&\quad+ \frac{c_0}{2}\|\partial_t\theta_{p}\|^2_{L^\infty(0,T;L^2)}+ \mu_p\|\mathbb{D}(\boldsymbol{\xi}_{p})\|^2_{L^\infty(0,T;L^2)}
+ \frac{\lambda_p}{2}\|\nabla\!\cdot\!\boldsymbol{\xi}_{p}\|^2_{L^2(0,T;L^2)}
\\
&\quad+ \mu_f\|\mathbb{D}(\partial_t\mathbf{u}_{f})\|^2_{L^2(0,T;L^2)} + \frac{a_f}{2}\|\nabla\partial_t\theta_{f}\|^2_{L^2(0,T;L^2)}
+ \frac{\mathbf{K}}{2}\|\nabla\partial_t p\|^2_{L^2(0,T;L^2)}
\\
&\quad+ \frac{\bm{a}_p}{2}\|\nabla\partial_t \theta_{p}\|^2_{L^2(0,T;L^2)}   \le (1+Te^{T}) \|\mathcal{C}_2\|_{L^2(0,T)}^2 + (1+Te^{T})\mathcal{C}_3^2,
\end{split}
\end{equation}
and
\begin{equation}
\begin{split}
\|p_f\|_{L^\infty(0,T;L^2)}
\le &\frac{1}{L}\Bigl(
C_{P,f}C_{K,f}\sqrt{2(1+Te^T)\big(\|\mathcal C_2\|_{L^2(0,T)}^2+\mathcal C_3^2\big)}
+ 2\mu_f\,\delta_1
+ C_{S,f}^2C_{K,f}^2\,\delta_1^2
\\
& + T_1T_3\sqrt{\frac{2(1+C_{P,p}^2)(1+e^T)}{K_{min}}\|\mathcal C_1\|_{L^2(0,T)}^2}
+ \frac{\alpha\sqrt{d}}{\sqrt{\operatorname{tr}(\mathbf{K})}} T_1^2C_{K,f}\,\delta_1
\\
& + \frac{\alpha\sqrt{d}}{\sqrt{\operatorname{tr}(\mathbf{K})}} T_1T_5
\sqrt{\frac{2(1+C_{P,p}^2)C_K^2(1+Te^T)}{\mu_p}\big(\|\mathcal C_2\|_{L^2(0,T)}^2+\mathcal C_3^2\big)}
+ \frac{G_r}{Re^2}C_{P,f}^2\,\delta_2\Bigr),
\end{split}
\label{eq:pf-Linfty}
\end{equation}
\end{theorem}
Here \(C_{S,f}, C_K, C_{P,f}, T_1, T_3, T_5\) denote the Sobolev embedding constant, Korn's constant, Poincaré constant, and trace constants, respectively.

\section{Conclusion}\label{sec:7}
In this work, we have established the well-posedness of a time-dependent
thermal-fluid-solid multiphysical coupled model that describes the interaction
between an incompressible free fluid and a thermo-poroelastic medium across a
fixed interface. Under a suitable small-data assumption, the existence of a weak
solution is proved via a Galerkin approximation scheme combined with uniform a
priori estimates and a compactness argument; uniqueness is established under
analogous conditions; and the fluid pressure is recovered through an inf-sup
condition, thereby completing the analysis of the full coupled system.

Several directions remain open for future investigation. First, the present
analysis relies on a small data assumption, which is used to control the
nonlinear convective terms and to guarantee the positivity of certain energy
functionals. Removing this restriction, or replacing it by a weaker
condition, is a challenging but important task. Second, the development and
analysis of fully discrete numerical schemes, together with optimal error
estimates, would turn the present well-posedness results into practical
tools for geothermal and biomechanical applications, where the coupling
between free-fluid flow, heat transfer, and thermo-poroelastic deformation
plays a central role.

\section*{Declarations}

\begin{itemize}
\item Conflict of interest\\
 The authors declare that they have no conflict of interest.
\end{itemize}

% References
\bibliographystyle{unsrt}

\bibliography{references}

\end{document}